\documentclass[12pt,a4paper,reqno]{amsart}
\usepackage{amssymb}
\usepackage{latexsym}
\usepackage{amsmath}
\usepackage{mathrsfs}
\usepackage{cite}

\usepackage{calc}

\newcommand{\scal}[2]{\langle #1,#2\rangle}
\newcommand{\rr}[1]{\mathbf R^{#1}}

\newcommand{\nm}[2]{\Vert #1\Vert _{#2}}

\newcommand{\sets}[2]{\{ \, #1\, ;\, #2\, \} }
\newcommand{\ep}{\varepsilon}
\newcommand{\fy}{\varphi}
\newcommand{\cdo}{\, \cdot \, }
\newcommand{\supp}{\operatorname{supp}}

\newcommand{\eabs}[1]{\langle #1\rangle}     

\newcommand{\vrum}{\vspace{0.1cm}}
\newcommand{\Char}{\operatorname{Char}}

\numberwithin{equation}{section}          

\newtheorem{thm}{Theorem}
\numberwithin{thm}{section}
\newtheorem{prop}[thm]{Proposition}

\newtheorem{lemma}[thm]{Lemma}
\newtheorem*{tom}{\rubrik}
\newcommand{\rubrik}{}

\theoremstyle{definition}

\newtheorem{defn}[thm]{Definition}
\newtheorem{example}[thm]{Example}

\theoremstyle{remark}
\newtheorem{rem}[thm]{Remark}

\author{Stevan Pilipovi\'c}

\address{Department of Mathematics and Informatics,
University of Novi Sad, Novi Sad, Serbia}

\email{stevan.pilipovic@dmi.uns.ac.rs}

\author{Nenad Teofanov}

\address{Department of Mathematics and Informatics,
University of Novi Sad, Novi Sad, Serbia}

\email{nenad.teofanov@dmi.uns.ac.rs}

\author{Joachim Toft}

\address{Department of Mathematics and Systems Engineering,
V{\"a}xj{\"o} University, Sweden}

\email{joachim.toft@vxu.se}

\title{\textbf {Micro-local analysis in Fourier Lebesgue and
modulation spaces. Part II}}

\begin{document}

\begin{abstract}
We consider different types of (local) products $f_1 f_2$ in Fourier
Lebesgue spaces. Furthermore, we prove the existence of such
products for other distributions satisfying appropriate wave-front
properties. We also consider  semi-linear equations of the form
$$
\qquad P(x,D)f = G(x,J_k f),
$$
with appropriate polynomials $P $ and $G$.
If the solution locally  belongs to appropriate weighted Fourier
Lebesgue space ${\mathscr F}L^q_{(\omega )} (\rr d)$ and $P$ is
non-characteristic at $(x_0,\xi_0),$ then we prove that
$(x_0,\xi_0)\not \in WF_{ {\mathscr F}L^q_{(\widetilde {\omega })} } (f)$,
where $\widetilde{\omega }$ depends on $\omega$, $P$ and $G$.
\end{abstract}

\maketitle

\section{Introduction}\label{sec0}

\par

In \cite{Hrm-nonlin}, H{\"o}rmander introduced wave-front sets with
respect to Sobolev spaces $ H^2_{s,loc}$ of Hilbert type, and used
such wave-front sets  to investigate regularity properties for
pseudo-differential operators, as well as solutions of semi-linear
equations. (We use the same notation for the usual function and
distribution spaces as in \cite{Ho1}.) H{\"o}rmander has proved
(\cite{Hrm-nonlin}) that for any appropriate pseudo-differential
operator $T$ with the set of characteristic points $\Char (T)$, and
for any distribution $f$, the wave-front set $WF _{H^2_s}(f)$ of the
distribution $f$ with respect to $H^2_{s}$ fulfills embeddings of
the form 
\begin{equation}\label{wavefrontrel1}
WF _{H^2_{s'}}(Tf) \subseteq WF _{H^2_s}(f) \subseteq WF
_{H^2_{s'}} (Tf)\bigcup \Char (T).
\end{equation}
Similar properties hold for the usual wave-front set (i.{\,}e.
wave-front set with respect to $C^\infty $ in \cite[Sections
8.1--8.3]{Ho1}), as well as for analytic wave-front sets. A relation
of the form \eqref{wavefrontrel1} is essential in the study of
regularity properties of solutions to partial differential
equations.

\par

Some of the results from  \cite{Hrm-nonlin} were extended in
different ways in \cite{PTT}, where general wave-front sets of
Fourier Lebesgue spaces were introduced. These (weighted) Fourier
Lebesgue spaces are denoted by $ \mathscr{F}L^{q} _{(\omega)} (\rr
d)$, where $ q \in [1,\infty]$ and $\omega$ is a weight function.
(See Section \ref{sec1} for strict definition.) In
particular, it is proved in  \cite{PTT} that \eqref{wavefrontrel1} holds for more
general type of wave-front sets. Furthermore, the class of
permitted pseudo-differential operators is larger comparing to
\cite{Hrm-nonlin}, and the set $\Char (T)$ can be replaced by a
smaller type of characteristic set, which better fits hypoelliptic
problems.

\par

In this paper we continue the study of \cite{PTT} concerning
local and microlocal properties of Fourier Lebesgue type spaces.
Especially, by using the framework of \cite{Hrm-nonlin}, we
establish multiplication properties for such wave-front sets, and
show how such results can be used to investigate the regularity of
solutions of semi-linear equations.

\par

We emphasize that our analysis is more subtle then the analysis in
\cite[Chapter 8]{Hrm-nonlin} since $ \mathscr{F}L^{q} _{(\omega)}
(\rr d)$, $ q\neq 2 $,  are not Hilbert spaces. In particular, the
analysis of cases $ 1\leq q \leq 2 $ and $ 2 < q $ is performed by
the use of different arguments, and we are forced to use extra
parameters in order to control the regularity. We refer to
\cite{PTT} and \cite{PTT3} for applications to the theory of
pseudo-differential operators, and here we apply our results to
the study of certain semi-linear equations.

\par


\par

We also give an interpretation of our results in the framework of
Feichtinger's modulation spaces (cf.
\cite{F1,Feichtinger3,Feichtinger4,Feichtinger5,Gro-book}).
Although those results are an immediate consequence of micro-local
analysis of $ \mathscr{F}L^{q} _{(\omega)} (\rr d)$, $ q \in
[1,\infty] $ and results from \cite{PTT}, they might be of an
independent interest in time-frequency analysis. For that reason,
microlocal results in context of modulation space theory are
collected in a separate section. Note that the wave front set with
respect to the modulation space $M^{p,q}_{(\omega)} (\rr d)$
coincides with the wave-front set with respect to $
\mathscr{F}L^{q} _{(\omega)} (\rr d) $ (see Section \ref{sec6} for
definitions). We refer to \cite{PTT} for the proof.

\par

The interest for such spaces in microlocal analysis has grown by
the recent work of several authors who have studied
pseudo-differential and Fourier integral operators within Fourier
Lebesgue spaces and their connection to modulation spaces in
different contexts (cf. \cite{CT1, CNR, Kasso, RSTT, Str, TCG1}).

\par

We end this introduction with a brief summary of the paper. In
Section \ref{sec1} we recall  basic facts of weighted Fourier
Lebesgue spaces and introduce the corresponding wave front sets.
We also discuss localized versions of  Fourier Lebesgue spaces and
relate the introduced wave front sets to H\"ormander's wave front
set (Proposition \ref{opadanje-za-sve-s} and Example
\ref{example}). Although the convolution is not the main objective
of the present paper and it will be studied in more details
elsewhere, in Section \ref{sec2} we prove continuity properties
and characterize the wave front set of the convolution (cf. Lemma
\ref{konvolucija}, Proposition \ref{convwavefront}).

\par

In Section \ref{sec3} we present three types of results on
products $f_1 f_2$, where $f_1$ and $f_2$ belong to appropriate
Fourier Lebesgue spaces (cf. Lemma \ref{inS'}, Proposition
\ref{p-q_estimates} and Theorem \ref{product}). Then, in Section
\ref{sec4}, we analyze microlocal properties  of products of
elements which locally belong to appropriate Fourier Lebesgue
spaces (cf. Theorems \ref{microlocal-functions} and
\ref{microlocal-functions-some-more}).

\par

In Section \ref{sec5} we use results from Sections \ref{sec3} and
\ref{sec4} to analyze microlocal regularity for a solution of
semi-linear equations of the form $ P(x,D)f = G(x,f, \dots,
f^{(k)})$, with appropriate assumptions on polynomials $P $ and $G$,
Theorems \ref{F} and \ref{semi-linear}. Furthermore, when $ q = 1,$
we are able to extend  the result of Theorem \ref{semi-linear} to
the case when $G$ is a real analytic function, Proposition
\ref{Fprop}. In Section \ref{sec6} we introduce modulation spaces
and restate some results in the terms of modulation spaces. In order
to be self-contained in Appendix A we recall necessary facts of a
class of  pseudo-differential operators which has been used in
Section \ref{sec5}.

\par

\section{Preliminaries}\label{sec1}

\par

In this section we make a review of notions and notation, and
discuss basic results.

\par

Assume that $x, \xi \in \rr d$. Then the scalar product of $x$ and
$\xi$ is denoted by $\langle x,\xi \rangle$, and we set $\langle x
\rangle = (1+ |x|^2)^{1/2}$. A conical neighborhood of a point $
(x_0, \xi_0) \in \rr d \times (\rr d\setminus 0)$ is a product $X
\times \Gamma $, where $X$ is an (open) neighborhood of $ x_0 $ in $
\rr d$ and $ \Gamma = \Gamma_{\xi_0}$ is an open cone in $ \rr d
\setminus 0 $ which contains $ \xi_0$. For $q \in [1,\infty]$ we let
$q'\in [1,\infty ]$ denote the conjugate exponent, i.{\,}e.
$1/q+1/q'=1$.

\par

The Fourier transform $\mathscr F$ is the linear and continuous
mapping on $\mathscr S'(\rr d)$ which takes the form
$$
(\mathscr Ff)(\xi )= \widehat f(\xi ) \equiv (2\pi )^{-d/2}\int
f(x)e^{-i\scal  x\xi }\, dx, \qquad \xi \in \rr d,
$$
when $f\in L^1(\rr d)$. We recall that $\mathscr F$ is a homeomorphism on
$\mathscr S'(\rr d)$ which restricts to a homeomorphism on $\mathscr
S(\rr d)$ and to a unitary operator on $L^2(\rr d)$.

\par

Assume that $\omega ,v\in L^\infty _{loc}(\rr d)$ are positive
(weight) functions. Then $\omega$ is called $v$-moderate if
\begin{equation}\label{moderate}
\omega (x+y) \leq C\omega (x)v(y)
\end{equation}
for some  constant $C$ which is independent of $x,y\in \rr
d$. If $v$ in \eqref{moderate} can be chosen as a polynomial, then
$\omega$ is called polynomially
moderated. We let $\mathscr P (\rr d)$ be the set of all polynomially
moderated functions on $\rr d$.

\par

Assume that $\omega \in \mathscr P(\rr {2d})$. Then the (weighted)
Fourier Lebesgue space $\mathscr FL^q_{(\omega )}(\rr d)$ is the
Banach space which consists of all $f\in \mathscr S'(\rr d)$ such that
\begin{equation}\label{FLnorm}
\nm f{\mathscr FL^{q}_{(\omega )}} = \nm f{\mathscr
FL^{q}_{(\omega ),x}} \equiv \nm {\widehat f\cdot \omega(x,\cdot )
}{L^q}<\infty .
\end{equation}

\par

\begin{rem}\label{unusual}
It might not seem to be natural to permit weights $\omega (x,\xi) $ in
\eqref{FLnorm} that are dependent on both $x$ and $\xi$, since
$\widehat f(\xi)$ only depends on $\xi$. We note that the fact
that $\omega$ is $v$-moderate for some $v\in \mathscr P(\rr {2d})$
implies that different choices of $x$ give rise to equivalent
norms. Therefore, the condition $\nm f{\mathscr FL^{q}_{(\omega
),x}}<\infty$ is independent of $x$.
\end{rem}

\par

Due to Remark \ref{unusual} we usually assume that the weights for the
Fourier Lebesgue spaces only depend on $ \xi $. Thus, with
$\omega_0 (\xi) = \omega (0, \xi) \in
\mathscr P(\rr d)$, we have
$$
f\in \mathscr FL^{q}_{(\omega)} (\rr d)= \mathscr FL^{q}_{(\omega _0)}
(\rr d)\quad \Longleftrightarrow \quad  \| f
\|_{\mathscr FL^{q}_{(\omega _0)}} \equiv \| \widehat f w_0 \|_{L^q} <
\infty .
$$

\par

The convention of indicating weight functions with parenthesis is
used also in other situations. For example, if $\omega \in
\mathscr P(\rr d)$, then $L^p_{(\omega )}(\rr d)$ is the set of
all measurable functions $f$ on $\rr d$ such that $f\omega \in
L^p(\rr d)$, i.{\,}e. such that $\nm f{L^p_{(\omega )}}\equiv \nm
{f\omega}{L^p}$ is finite.

\par

Let $X $ be an open set in $ \rr d$. Then the \emph{local} Fourier
Lebesgue space $\mathscr FL^q_{(\omega ),loc}(X)$ consists of all
$f\in \mathscr D'(X)$ such that $\fy f\in \mathscr FL^q_{(\omega
)}(\rr d)$ for each $\fy \in C_0^\infty (X)$. The topology in
$\mathscr FL^q_{(\omega ),loc}(X)$  is defined by the family of
seminorms $f\mapsto  \nm {\fy f}{\mathscr FL^q_{(\omega )}}$,
where $\fy \in C_0^\infty (X)$.

\par

We note that if $\omega \in \mathscr P(\rr d)$, then
\begin{equation}\label{FLqFLqlocemb}
\mathscr FL^q_{(\omega )}(\rr d)\subseteq \mathscr FL^q_{(\omega
),loc}(X).
\end{equation}
In fact, if $f\in \mathscr FL^q_{(\omega )}(\rr d)$,  $\fy  \in
C_0^\infty (X)$ and if $v\in \mathscr P(\rr d)$ is chosen such that
$\omega $ is $v$ moderate, then Young's inequality gives
\begin{multline*}
\nm {\fy f}{\mathscr FL^q_{(\omega)}} = \nm {\mathscr F(\fy f)\,
\omega}{L^q} =(2\pi )^{-d/2} \nm {(\widehat \fy *\widehat f \, )\,
\omega}{L^q}
\\[1ex]
\le (2\pi )^{-d/2} \nm {|\widehat \fy \, v|* |\widehat f\, \omega
|}{L^q} \le C\nm {\widehat f\, \omega}{L^q} =C\nm {f}{\mathscr
FL^q_{(\omega)}},
\end{multline*}
where $C=(2\pi )^{-d/2}\nm {\widehat \fy \, v}{L^1}<\infty$. This
proves \eqref{FLqFLqlocemb}.

\par

Next we show that $\mathscr FL^q_{(\omega ),loc} (X)$ increases
with $q$ and decreases with $\omega$, i.{\,}e.
\begin{equation}\label{incrFLloc}
\mathscr FL^{q_1}_{(\omega _1),loc}(X)\subseteq \mathscr
FL^{q_2}_{(\omega _2),loc}(X), \;\; \text{when} \;\; q_1\le q_2 \;\;
\text{and} \;\; \omega _2\le C\omega _1.
\end{equation}

\par

The decrease with respect to $\omega = \omega_0$ is
straightforward. It remains to show the increase with respect to
$q$. Assume, without any loss of generality, that $f\in \mathscr
E'(X)$ and $\fy \in C_0^\infty (\rr d)$ are such that $\fy \equiv 1$
in $\supp f$, and choose $p\in [1,\infty ]$ such that
$1/q_1+1/p=1/q_2+1$. Then, for a $v$-moderate weight $\omega$, it
follows from Young's inequality that
\begin{multline*}
\nm {f}{\mathscr FL^{q_2}_{(\omega_0)}} \le C\nm {(\widehat \fy
*\widehat f\, ) \omega_0}{L^{q_2}}
\\[1ex]
\le C\nm {|\widehat \fy v|*|\widehat f \omega|}{L^{q_2}}
\le C\nm {\widehat \fy v}{L^p}\nm {\widehat f\omega}{L^{q_1}},
\end{multline*}
for some positive  constant $C$. Since $\widehat \fy \in \mathscr
S (\rr d)$, it follows that
\begin{equation*}
\nm {f}{\mathscr FL^{q_2}_{(\omega_0)}} \le C\nm {f}{\mathscr
FL^{q_1}_{(\omega_0)}},\quad \text{when}\quad q_1\le q_2,\quad f\in
\mathscr E'(X),
\end{equation*}
and the assertion follows.

\par

When $\omega_0 (\xi )= \eabs \xi ^s$
it is convenient to set
$$
\mathscr FL^q_s (\rr d) = \mathscr FL^q_{(\omega_0 )} (\rr d) ,\quad
\mathscr FL^q_{s,loc} (X) = \mathscr FL^q_{(\omega_0 ),loc} (X).
$$
Furthermore, if $\omega (x)= \eabs x^s $, $ s \in {\bf R}$,
then we use the notation $L^p _s (\rr d) $ instead of $L^p_{(\omega
)} (\rr d) $. We also omit the indices of weights when $\omega
\equiv 1$, i.{\,}e. we set
\begin{equation*}
\mathscr FL^{q}= \mathscr FL^{q}_{(\omega )},
\quad \mathscr FL^{q}_{loc} = \mathscr FL^{q}_{(\omega ),loc} ,\quad
L^{p}  =L^{p}_{(\omega )},
\quad \text{when}\quad
\omega \equiv 1.
\end{equation*}

\par

Recall that if $ f \in \mathscr{D}' (\rr d) $ then the wave front
set of $f$, $WF (f)$ is the set of points $ (x_0,\xi_0) \in \rr d
\times (\rr d \setminus 0)$ such that for no conical neighborhoods $
X \times \Gamma_{\xi_0} $ of $(x_0,\xi _0)$, no $ \chi \in
C^{\infty} _0 (X)$ with $ \chi (x_0)\neq 0,$ and no constants $C_N$
depending on $N$ only, we have
$$
| \mathscr F(\chi f) (\xi) | \leq C_N \langle \xi \rangle
^{-N/2},\quad \xi \in \Gamma_{\xi_0},\quad \text{and}\quad N \in
{\bf N}.
$$

\par

Next we define  wave front set with respect to Fourier Lebesgue spaces.

\par

Assume that $\omega \in \mathscr P(\rr {2d})$, $\Gamma
\subseteq \rr d\setminus 0$ is an open cone and $q\in [1,\infty ]$
are fixed. For any $f\in \mathscr S'(\rr d)$, let
\begin{equation}
|f|_{\mathscr FL^{q,\Gamma}_{(\omega ), x}}
\equiv
\Big ( \int _{\Gamma} |\widehat f(\xi )\omega (x,\xi )|^{q}\, d\xi
\Big )^{1/q}\label{skoff}
\end{equation}
(with obvious interpretation when $q=\infty$). We note that $|\cdo
|_{\mathscr FL^{q,\Gamma}_{(\omega ), x}}$ defines a semi-norm on
$\mathscr S'(\rr d)$ which might attain the value $+\infty$. Since
$\omega $ is $v$-moderate for some $v \in \mathscr P(\rr {2d})$,
it follows that different $x \in \rr d$ gives rise to equivalent
semi-norms $ |f|_{\mathscr FL^{q,\Gamma}_{(\omega ), x}}$.
Furthermore, if $\Gamma =\rr d\setminus 0$, $f\in \mathscr
FL^{q}_{(\omega )}(\rr d)$ and $q<\infty$, then $|f|_{\mathscr
FL^{q,\Gamma}_{(\omega ), x}}$ agrees with the Fourier Lebesgue
norm $\nm f{\mathscr FL^{q}_{(\omega ), x}}$ of $f$.

\par

We let $ \Theta_{\mathscr FL^{q} _{(\omega )}}
(f)$ be the set of all $ \xi \in \rr d\setminus 0 $ such that
$|f|_{\mathscr FL^{q,\Gamma}_{(\omega ), x}} < \infty$, for some $
\Gamma = \Gamma_{\xi}$.  We also let $\Sigma_{\mathscr
FL^{q} _{(\omega )}} (f)$ be the complement of $ \Theta_{\mathscr
FL^{q} _{(\omega )}} (f)$ in $\rr d\setminus 0 $. Then
$\Theta_{ \mathscr FL^{q} _{(\omega )}} (f)$ and $\Sigma_{\mathscr
FL^{q} _{(\omega )}} (f)$ are open respectively
closed subsets in $\rr d\setminus 0$, which are independent of
the choice of $ x \in \rr d$ in \eqref{skoff}. We have now the
following result, see \cite{PTT}.

\par

\begin{prop}\label{theta-sigma-propertiesAA}
Assume that  $q \in [1,\infty]$, $ \chi \in  \mathscr S(\rr d)$, and
that $\omega \in \mathscr P(\rr {2d})$. Also assume that $f\in
\mathscr E'(\rr d)$. Then
\begin{equation}\label{chi-subsetAA}
\Sigma_{\mathscr FL^{q}_{(\omega )} } (\chi f) \subseteq
\Sigma_{\mathscr  FL^{q}_{(\omega )} } (f).
\end{equation}
\end{prop}

\par

\begin{defn}\label{wave-frontdef}
Assume that  $q\in [1,\infty ]$, $X\subseteq \rr d$ is open,
$f\in \mathscr D'(X)$, $\omega _0\in \mathscr P(\rr d)$ and
$\omega \in \mathscr P(\rr {2d})$ are such that $\omega _0(\xi
)=\omega (y_0,\xi )$ for some $y_0\in \rr d$. The wave-front set
$$
WF_{\mathscr FL^q_{(\omega _0)}}(f)\equiv WF_{\mathscr FL^q_{(\omega
)}}(f)
$$
with respect to $\mathscr F L^{q} _{(\omega_0 )}(\rr d)=\mathscr F
L^{q} _{(\omega )}(\rr d)$ consists of all pairs $(x_0,\xi_0)$ in
$X\times (\rr d \setminus 0)$ such that
$$
\xi _0 \in \Sigma_{\mathscr FL^{q} _{(\omega _0)}} (\chi
f) = \Sigma_{\mathscr FL^{q} _{(\omega )}} (\chi f)
$$
holds for each $\chi \in C_0^\infty (X)$ such that $\chi (x_0)\neq
0$.
\end{defn}

\par

We note that $WF_{\mathscr FL^q_{(\omega _0)}}(f)$ in Definition
\ref{wave-frontdef} is a closed set in $X\times (\rr d\setminus
0)$, which is independent of the choice of $y_0\in \rr d$
and $p\in [1,\infty ]$, since its complement is open in $X\times
(\rr d\setminus 0)$.

\par

The following proposition shows that the wave-front set
$WF_{\mathscr FL^q_{(\omega )} }(f)$ decreases with respect to the
parameter $q$ and increases with respect to the weight function
$\omega$, when $f\in \mathscr D'(X)$ is fixed. We refer to
\cite[Proposition 2.3]{PTT} for the proof.

\par

\begin{prop}\label{wavefrontprop1}
Assume that $X$ is an open subset of $\rr d$, $f\in \mathscr D'(X)$,
$q_j\in [1,\infty ]$ and $\omega _j\in \mathscr P(\rr {2d})$ for
$j=1,2$ satisfy
\begin{equation}\label{qandomega}
q_1\le q_2,\quad \text{and}\quad \omega _2(x,\xi )\le C\omega
_1(x,\xi ),
\end{equation}
for some postivie constant $C$ which is independent of $x,\xi \in \rr d$.
Then
$$
WF_{\mathscr FL^{q_2}_{(\omega _2)} }(f)\subseteq
WF_{\mathscr FL^{q_1}_{(\omega _1)} }(f).
$$
\end{prop}

\par

We remark that several properties for wave-front sets of Fourier
Lebesgue types can be found in \cite{PTT}. For example, it follows
from \cite[Theorem 3.1]{PTT} that
\begin{equation}\label{WFderiv}
WF _{\mathscr FL^{q}_{(\omega \vartheta )} }(\partial
_jf)\subseteq WF _{\mathscr FL^{q}_{(\omega )} }(f),\qquad
\vartheta (\xi )=(1+|\xi |^2)^{1/2}.
\end{equation}

\par

The following proposition describes the relation between the classical
wave-front set, analytic wave-front sets and wave-front sets of
Fourier Lebesgue types. Here recall that $WF_A(f)$ denotes the
analytic wave-front set (cf. \cite[Section 8.4]{Ho1}) of a
distribution $f$.

\par

\begin{prop} \label{opadanje-za-sve-s}
Let $q\in [1,\infty ]$ $X\subseteq \rr d$ be open, $ f\in \mathscr{D}'(X)$
and $ (x_0,\xi_0) \in X \times (\rr d \setminus 0).$ Then the
following conditions are equivalent:
\begin{enumerate}
\item[(1)] $ (x_0,\xi_0) \not \in WF_{\mathscr{F}L^{q} _{(\omega_0)}}
(f) $;

\vrum

\item[(2)] there exists $ g \in  \mathscr{F}L^{q} _{(\omega_0)} (\rr
d)$ ($ g \in  \mathscr{F}L^{q} _{(\omega_0),loc} (X)$) such that
$(x_0,\xi_0) \not \in WF (f - g) $;

\vrum

\item[(3)] there exists $ g \in  \mathscr{F}L^{q} _{(\omega_0)} (\rr
d)$ ($ g \in  \mathscr{F}L^{q} _{(\omega_0),loc} (X)$) such that $
(x_0,\xi_0) \not \in WF _A(f - g) $.
\end{enumerate}
\end{prop}

\par

For the proof we need the following result on multiplications of
elements in Fourier Lebesgue spaces. Here the involved exponents
should satisfy
\begin{align}
\frac 1{q_1}+\frac 1{q_2} &= 1+\frac 1{q}\label{youngcond3q}
\intertext{or}
\frac 1{q_1}+\frac 1{q_2} &\ge 1+\frac
1{q}.\tag*{(\ref{youngcond3q})$'$}
\end{align}

\par

\begin{lemma}\label{inS'B}
Assume that
$X\subseteq \rr d$ is open, $q,q_1,q_2\in [1,\infty ]$ and
that $\omega ,\omega _1,\omega _2\in \mathscr P(\rr d)$ satisfy
\begin{equation}\label{youngcondandmoreAA}
\omega (\xi _1+\xi _2)\le
C\omega _1(\xi _1)\omega _2(\xi _2),
\end{equation}
for some  constant $C$ which is independent of $\xi _1,\xi
_2\in \rr d$. Then the following is true:
\begin{enumerate}
\item if \eqref{youngcond3q} holds, then the map $(f_1,f_2)\mapsto f_1
f_2$ from $\mathscr S(\rr d)\times \mathscr S(\rr d)$ to $\mathscr
S(\rr d)$ extends uniquely to a continuous mapping from $\mathscr
FL^{q_1}_{(\omega _1)}(\rr d)\times \mathscr FL^{q_2}_{(\omega
_2)}(\rr d)$ to $\mathscr FL^{q}_{(\omega )}(\rr d)$;

\vrum

\item if \eqref{youngcond3q}$'$ holds, then  the map $(f_1,f_2)\mapsto
f_1 f_2$ from $C_0^\infty(X)\times C_0^\infty (X)$ to $C_0^\infty
(X)$ extends uniquely to a continuous mapping from $\mathscr
FL^{q_1}_{(\omega _1),loc}(X)\times \mathscr FL^{q_2}_{(\omega
_2),loc}(X)$ to $\mathscr FL^{q}_{(\omega ),loc}(X)$.
\end{enumerate}
\end{lemma}

\par

\begin{proof}
(1) We must have that $q_1<\infty$ or $q_2<\infty$ in view of
\eqref{youngcond3q}. Therefore assume that
$q_2<\infty$, and let $f_1\in \mathscr FL^{q_1}_{(\omega _1)} (\rr d)$
and $f_2\in \mathscr S (\rr d)$.

\par

By \eqref{youngcondandmoreAA} it follows that
\begin{equation*}
|(\widehat f_1*\widehat f_2)(\xi )\omega (\xi )| \le C (|\widehat
f_1 \omega _1|*|\widehat f_2\omega _2|)(\xi ),
\end{equation*}
for some  constant $C$. Hence, by Young's inequality we get
\begin{equation*}
\nm {f_1f_2}{\mathscr FL^q_{(\omega )}}\le C\nm {\, |\widehat f_1
\omega _1|*|\widehat f_2\omega_2|\, }{L^q}\le C \nm {\widehat f_1
\omega _1 }{L^{q_1}}\nm {\widehat f_2 \omega _2}{L^{q_2}}.
\end{equation*}
The assertion (1) now follows from this estimate
and the fact that $\mathscr S (\rr d)$ is dense in $\mathscr
FL^{q_2}_{(\omega _2)} (\rr d)$, since $q_2<\infty$.

\par

The assertion (2) is an immediate consequence of \eqref{FLqFLqlocemb},
\eqref{incrFLloc} and (1). The proof is complete.
\end{proof}

\par

\begin{proof}[Proof of Proposition \ref {opadanje-za-sve-s}.]
The idea for the proof is the same as for the proof of
\cite[Proposition 8.2.6]{Hrm-nonlin} (see also  \cite{PTT3}). For the
sake of completeness, and since similar arguments will be used later
on, we give a proof here.

\par

Since wave-front sets are locally defined, we may assume $g\in
\mathscr FL^{q}_{(\omega _0),loc}(X)$ in (2) and (3). Assume that
(2) holds. Then it follows that for some open subset $X_0$ of $X$
and for some $\Gamma =\Gamma _{\xi _0}$ we have
\begin{equation}\label{locfourest1}
|\mathscr F(\chi (f-g))(\xi )|\le C_{N,\chi }\eabs \xi ^{-N},\quad
N=1,2,\dots ,\quad \xi \in \Gamma ,
\end{equation}
when $\chi \in C_0^\infty (X_0)$. In particular it follows that
$|\chi (f-g)|_{\mathscr FL^{q,\Gamma}_{(\omega _0)}}$ is finite.
Since $g \in \mathscr{F}L^{q} _{(\omega_0)} (\rr d)$, it follows
that $|\chi g|_{\mathscr FL^{q,\Gamma}_{(\omega _0)}} $ is finite,
by Lemma \ref{inS'B} (1). This implies that $|\chi f|_{\mathscr
FL^{q,\Gamma}_{(\omega _0)}} < \infty$ for each $\chi \in C_0^\infty
(X_0)$, and (1) holds.

\par

Conversely, if $  (x_0,\xi_0) \notin WF_{\mathscr{F}L^{q}
_{(\omega_0)}} (f)$, then there exists a conical neighborhood $X_0
\times \Gamma$ of $(x_0,\xi _0)$ such that
\begin{equation}\label{varphi-u}
|\chi f|_{\mathscr FL^{q,\Gamma}_{(\omega _0)}}<\infty ,\quad
\text{when}\quad \chi \in C_0^\infty (X_0),
\end{equation}
in view of Proposition \ref{theta-sigma-propertiesAA}.

\par

Now let $\chi _1,\chi \in C_0^\infty (X_0)$ be chosen such that
$\chi (x_0)\neq 0$ and $\chi _1=1$ in the support of $\chi $,
and let $g$ be defined by the formula
$$
\widehat g (\xi) = \left \{
\begin{array}{lll}
\mathscr F(\chi _1f) (\xi), & \mbox{if} & \xi \in \Gamma
\\
0, & \mbox{if} & \xi \not\in \Gamma .
\end{array}
\right .
$$
Then $g\in \mathscr FL^q_{(\omega _0)}(\rr d)$ and
\begin{equation}\label{stiligt}
|\mathscr F(\chi _1f-g)(\xi )|\le C_N\eabs \xi ^{-N},\quad N=0,1,2,\dots ,
\end{equation}
when $\xi \in \Gamma$, since the left-hand side in \eqref{stiligt} is
identically equal to zero.

\par

By \cite[Lemma 8.1.1]{Ho1} and its proof, it follows that
\eqref{stiligt} holds after $\chi f-g$ is replaced by $\chi _1(\chi
f-g)$, and $\Gamma$ is replaced by a smaller conical neighborhood of
$\xi _0$, if necessary. Since $\chi _1\chi =\chi _1$, it follows
that \eqref{locfourest1} holds, which means that $(x_0,\xi _0)\notin
WF(f-g)$. This proves the equivalence between (1) and (2).

\par

Since $WF(f)\subseteq WF_A(f)$ for each distribution $f$, it
follows that (2) holds when (3) is fulfilled. Assume instead that
(2) holds. Then for some $h\in C^\infty (X)$ we have $(x_0,\xi
_0)\notin WF_A(f-g-h)$, in view of the remark before
\cite[Corollary 8.4.16]{Ho1}. Since $C^\infty (X)\subseteq
\mathscr FL^q_{(\omega _0),loc}(X)$, it follows that $g_1=g+h\in
\mathscr FL^q_{(\omega _0),loc}(X)$. Hence (3) holds, and the
result follows.
\end{proof}

\par

Next we recall wave-front sets of superior
type, which, together with wave-front sets of inferior
type, were considered in \cite{PTT}. Let $q \in [1,\infty]$. For each $s_0\in \mathbf
R\cup \{ \infty \}$ and $f\in \mathscr E'(\rr d)$, we let $ \Theta_{\mathscr
FL^{q} _{s_0}}^{\sup} (f)$ be the set of all $ \xi \in \rr d\setminus
0 $ such that for some $ \Gamma = \Gamma_{\xi}$ and each $s<s_0$ we
have $|f|_{\mathscr FL^{q,\Gamma}_{s}} < \infty$. We also let $ \Sigma
_{\mathscr FL^{q} _{s_0}}^{\sup} (f)$ be the complement in $\rr
d\times (\rr d\setminus 0)$ of $ \Theta _{\mathscr
FL^{q}_{s_0}}^{\sup} (f)$.

\par

\begin{defn} \label{glob}
Assume that $s_0 \in \mathbf R\cup  \{ \infty \}$, and let $ f\in
\mathscr{D}'(X)$. Then the wave-front set of superior type,
$WF^{\operatorname{sup}}_{\mathscr FL^q_{s_0}}(f)$ of $f$, consists of
all pairs $(x_0,\xi _0)$ in $X\times (\rr d\setminus 0)$ such that for no
$\chi \in C_0^\infty (X)$ with $\chi (x_0)\neq 0$, and no $\Gamma
=\Gamma _{\xi _0}$ it holds $|\chi f|_{\mathscr FL^{q,\Gamma } _{s}}
<\infty$, for each $s<s_0$.
\end{defn}

\par

It follows that $WF^{\operatorname{sup}}_{\mathscr FL^q_{s_0}}(f)$ consists
of all pairs $(x_0,\xi _0)$ in $X\times (\rr d\setminus 0)$ such
that $\xi _0 \in \Sigma_{\mathscr FL^{q} _{s_0}} ^{\sup} (\chi f)$
holds for each $\chi \in C_0^\infty (X)$ such that $\chi (x_0)\neq
0$. The following
result follows from Remark 4.2 in \cite{PTT}.

\par

\begin{prop} \label{istiset}
Let  $ f\in \mathscr{D}'(X).$ Then
$ WF^{\operatorname{sup}}_{\mathscr FL^q_\infty}(f) =  WF(f)$.
\end{prop}

\par

In the definition of $ WF^{\operatorname{sup}}_{\mathscr FL^q_\infty}(f) $
the corresponding cone $\Gamma $ is fixed and independent on $s.$
However, one may consider the following situation.

\par

\begin{prop} \label{prop-loc}
Let $ f \in \mathscr{D}'(X),$ $ X $ be  open  in $ \rr d$ and $
(x_0,\xi_0) \in X \times (\rr d \setminus 0)$.
Then the following conditions are equivalent:
\begin{enumerate}
\item[(1)] $(x_0,\xi _0)\notin \cup_{s\geq 0} WF_{\mathscr{F}L^{q}
_{s}} (f)$;

\vrum

\item[(2)]  for every $s\geq 0$
there exists an open subset $ X_s \subset X$ and
$g_s \in \mathscr{F}L^{q} _{s, loc} (X_s) $ such that
$ (x_0,\xi_0) \not \in WF (f-g_s)$;

\vrum

\item[(3)] for every $s\geq 0$, there exists $\chi _s \in C^\infty _0
(X)$, $ \chi _s \equiv 1$ in a neighborhood of $ x_0$, and a cone $
\Gamma _s = \Gamma_{s, \xi_0}$ such that $|\chi _sf|_{\mathscr
FL^{q,\Gamma _s}_{s}}<\infty$.
\end{enumerate}
\end{prop}

\par

Proposition \ref{prop-loc} follows immediately from Definition
\ref{wave-frontdef}, Proposition \ref{wavefrontprop1} and the
proof of Proposition \ref{opadanje-za-sve-s}.

\par

The next proposition follows immediately if we replace $\Gamma_s$
with $\rr d$ in Proposition \ref{prop-loc} (3), and use
compactness arguments.

\par

\begin{prop}
Let $ f \in \mathscr{D}'(X)$ and assume that the conditions
{\rm{(1)--(3)}} in Proposition \ref{prop-loc} hold for some $ x_0
\in X $ and every $\xi_0 \in \rr d \setminus 0$. Then for every
$s\geq 0$ there exists $\chi_s \in C^\infty _0 (X),$ $ \chi _s
(x_0) \neq 0$ such that $\chi _sf \in \mathscr FL^q_s$.
\end{prop}

\par

By the above definitions it follows that
\begin{equation}\label{cupWFemb}
\cup_{s\geq 0} WF_{\mathscr{F}L^{q} _{s}} (f) \subseteq   WF (f).
\end{equation}
The next example shows that the embedding in \eqref{cupWFemb}
might be strict.

\par

\begin{example} \label{example}
Let $q\in [1,\infty]$ and $ ( [a_n, b_n] )_{n \in {\mathbf N}}$ be a
sequence of disjoint intervals so that $ a_n \searrow 0 $ and $
b_n < a_{n+1}$, $n \in \mathbf N$. Also let $ f_n \in \mathscr{F}L^{q}
_{n+2} (\mathbf R)\setminus \mathscr{F}L^{q} _{n+3}(\mathbf R)$ be
such that $ \supp f_n \subseteq [a_n, b_n]$ and let $ M_n \equiv
\sup_{x \in [a_n, b_n]} | f_n (x) |$, $ n\in \mathbf N$. Then
$M_n<\infty$, since $ \mathscr{F}L^{q} _{n+2}(\mathbf R) \subseteq
\mathscr{F}L^{1}(\mathbf R)\subseteq L^\infty (\mathbf R)$, when
$n\ge 0$. In particular,
$$
f \equiv \sum \limits_{n=0}^\infty \frac{1}{n^2 M_n} f_n
$$
is well-defined. By the construction we have that $(0, \xi_0) \not\in
\cup_{s\geq 0} WF_{\mathscr{F}L^{q} _{s}} (f)$ for any fixed $ \xi_0
\in \mathbf R \setminus 0$. On the other hand, $(0, \xi_0) \in  WF
(f)$, when $ \xi_0 \in \mathbf R \setminus 0$.
\end{example}

\par

\section{Convolution properties for wave-front sets in Fourier
Lebesgue spaces} \label{sec2}

\par

In this section we prove that convolution properties, valid for
standard wave-front sets of H{\"o}rmander type, also hold for the
wave-front sets of Fourier Lebesgue types. We first
consider convolutions for distributions in Fourier Lebesgue spaces.

\par

\begin{lemma} \label{konvolucija}
Assume that $q,q_1,q_2\in [1,\infty ]$ and $ \omega, \omega_1,
\omega_2 \in \mathscr P (\rr d)$ satisfy
\begin{equation}\label{qomegarel1}
\frac 1{q_1}+\frac 1{q_2} = \frac 1{q}\quad \text{and}\quad  \omega
(\xi) \leq  C \omega_1 (\xi) \omega_2 (\xi),
\end{equation}
for some $ C> 0 $ which is independent of $\xi \in \rr d$. Then the
convolution map $(f_1,f_2)\mapsto f_1*f_2$ from $\mathscr S(\rr
d)\times \mathscr S(\rr d)$ to $ \mathscr S (\rr d)$ extends to a
continuous mapping from $\mathscr FL^{q_1}_{(\omega_1)}(\rr d)\times
\mathscr FL^{q_2}_{(\omega_2)}(\rr d)$ to $\mathscr
FL^{q}_{(\omega)}(\rr d)$. This extension is unique if $q_1<\infty$ or
$q_2<\infty$.
\end{lemma}

\par

\begin{proof}
First we assume that $q_2<\infty$, and we let $f\in \mathscr
FL^{q_1}_{(\omega_1)}(\rr d)$ and $f_2\in \mathscr S(\rr d)$. Then
$q<\infty$, and $f_1*f_2$ is well-defined as an element in $\mathscr
S'(\rr d)$. By \eqref{qomegarel1} and H{\"o}lder's inequality we get
\begin{equation*}
\nm {f_1*f_2}{\mathscr F L^q_{(\omega)}}
=
(2\pi )^{d/2}\nm {\widehat f_1\widehat f_2\omega}{L^q}
\le
C \nm {\widehat f_1 \omega_1 }{L^{q_1}} \nm
{\widehat f_2 \omega_2 }{L^{q_2}},
\end{equation*}
for some positive constant $C$.
This gives
\begin{equation}\label{convest1}
\nm {f_1*f_2}{ \mathscr FL^q _{(\omega)} } \le (2\pi )^{d/2} C
\nm {f_1}{\mathscr FL^{q_1} _{(\omega_1)} } \nm {f_2}{ \mathscr
FL^{q_2} _{(\omega_2)}}.
\end{equation}
The result now follows in this case from the fact that $\mathscr S
$ is dense in $\mathscr FL^{q_2} _{s_2}$, and by similar
arguments, the result follows when $q_1<\infty$.

\par

It remains to consider the case $q_1=q_2=\infty$. We note that
\begin{equation}\label{convreform}
f_1*f_2=(2\pi )^{d/2}\mathscr F(\widehat f_1\cdot \widehat f_2)
\end{equation}
holds when $f_1,f_2\in \mathscr S(\rr d)$. The asserted extension
follows if we prove that the right-hand side of \eqref{convreform}
exists for $f_j\in \mathscr FL^{\infty}_{(\omega_j)}(\rr d)$, $ j
=1,2,$ and defines an element in $\mathscr FL^{\infty}_{(\omega)}(\rr
d)$.

\par

By the assumption we have  that $\widehat f_1$ and
$\widehat f_2$ are measurable and essentially bounded by some
polynomials. Hence the product $\widehat f_1\widehat f_2$ is
well-defined as a measurable function and essentially bounded by an
appropriate polynomial. By \eqref{convreform} it follows that
$f_1*f_2$ is well-defined as an element in $\mathscr S'(\rr
d)$. Furthermore,
\begin{multline*}
\nm {f_1*f_2}{\mathscr FL^\infty _{(\omega )}} = (2\pi )^{d/2}\nm
{\widehat f_1\, \widehat f_2\omega}{L^\infty }
\le (2\pi )^{d/2} C \nm {(\widehat f_1\omega_1 )\, (\widehat
f_2\omega_2 )}{L^\infty }
\\[1ex]
\le (2\pi )^{d/2} C \nm {\widehat f_1\omega_1}{L^\infty}\nm
{\widehat f_2\omega_2}{L^\infty} = C \nm {f_1}{\mathscr
FL^\infty _{(\omega_1)}}\nm {f_2}{\mathscr FL^\infty _{(\omega_2)}},
\end{multline*}
and the result follows.
\end{proof}

\par

\begin{prop} \label{convwavefront}
Assume that $q,q_1,q_2\in [1,\infty ]$ and $ \omega, \omega_1,
\omega_2 \in \mathscr P (\rr d)$ satisfy \eqref{qomegarel1}, and that
$f_1\in \mathscr FL^{q_1}_{(\omega _1),loc}(\rr d)$ and $f_2\in \mathscr
D'(\rr d)$ are such that $f_1$ or $f_2$ has compact support. Then
$$
WF_{\mathscr FL^{q}_{(\omega )}}(f_1*f_2) \subseteq \sets {(x+y,\xi
)}{x\in \supp f_1\ \text{and}\ (y,\xi )\in WF_{\mathscr
FL^{q_2}_{(\omega _2)}}(f_2)}.
$$
\end{prop}

\par

For the proof we let $B_r(x_0)$ be the open ball in $\rr d$ with
radius $r>0$ and center at $x_0\in \rr d$.

\par

\begin{proof}
From the local property of the wave-front sets (cf. Proposition 2.1 in
\cite{PTT}), and the fact that one of $f_1$ and $f_2$ has compact
support, it follows from \eqref{FLqFLqlocemb} that we may assume that
$f_1\in \mathscr E' (\rr d) \cap \mathscr FL^{q_1}_{(\omega _1)} (\rr
d)$ and $f_2\in \mathscr E' (\rr d)$.

\par

Let $(x_0,\xi _0)$ be chosen such that $(y,\xi _0)\notin WF_{\mathscr
FL^{q_2}_{(\omega _2)}}(f_2)$ when $x\in \supp f_1$ and $x_0=x+y$, and
let $F(x,t)=f_2(x-t)f_1(t)$. Since $f_1$ and $f_2$ have compact
support, it follows by the definition that for some $\Gamma =\Gamma
_{\xi _0}$,
\begin{gather*}
y_1,\dots ,y_n\in \rr d,\qquad r>0,\qquad r_0>0,
\\[1ex]
\chi \in C_0^\infty (B_r(0)),\qquad \chi _0 \in C_0^\infty
(B_{r_0}(x_0)),\qquad \text{and}\quad \chi _1\in C_0^\infty (\rr d),
\end{gather*}
with $\chi _0=1$ in a neighborhood of $x_0$, we have
\begin{align}
|\chi (\cdo -y_j)f_2|_{\mathscr FL^{q_2,\Gamma}_{(\omega _2)}}
&<\infty ,\qquad 1\le j\le n\label{notinfinity}
\intertext{and}
\sum _{j=1}^n\chi (x-t-y_j)\chi _1(t)& \equiv 1 \notag
\end{align}
when $x\in \supp \chi _0$ and $t$ belongs to
$$
\sets {t\in \rr d}{(x,t)\in \supp F \ \text{for some}\ x\in \supp \chi
_0}.
$$

\par

By H{\"o}lder's inequality we get
\begin{multline*}
|\chi _0(f_1*f_2)|_{\mathscr FL^{q,\Gamma}_{(\omega )}} \le C \sum
_{j=1}^n |((\chi _1f_1)*(\chi (\cdo -y_j)f_2)|_{\mathscr
FL^{q,\Gamma}_{(\omega )}}
\\[1ex]
\le C\sum _{j=1}^n \nm {(\mathscr F(\chi _1f_1)\omega _1)\, (\mathscr
F(\chi (\cdo -y_j)f_2)\omega _2)}{L^q(\Gamma)}
\\[1ex]
\le C \sum _{j=1}^n \nm {\mathscr F(\chi _1f_1)\omega
_1}{L^{q_1}(\Gamma)} \nm {\mathscr F(\chi (\cdo -y_j)f_2)\omega
_2}{L^{q_2}(\Gamma)},
\end{multline*}
for some positive constant $C$.
Since the right-hand side is finite in the view of \eqref{notinfinity},
it follows that $(x_0,\xi _0)\notin WF_{\mathscr FL^{q}_{(\omega )}}
(f_1*f_2)$, and the result follows.
\end{proof}

\par

\begin{rem}
Assume that $f_1\in \mathscr D'(\rr d)$ and $f_2\in \mathscr E'(\rr d)$. Proposition \ref{convwavefront} in combination with the fact that
$$
\cup _{\omega \in \mathscr P}\mathscr FL^q_{(\omega ),loc} (\rr d)=\mathscr D'(\rr d)
$$
and Proposition \ref{istiset} can be used to prove
$$
WF(f_1*f_2)\subseteq \sets {(x+y,\xi )}{(x,\xi )\in WF(f_1)\ \text{and}\ (y,\xi )\in WF(f_2)}.
$$
\end{rem}

\par

\section{Multiplication in Fourer-Lebesgue
spaces}\label{sec3}

\par

In this section we discuss the problem of multiplication  in
$\mathscr F L^q_{s} (\rr d)$ and $\mathscr F L^q_{s,loc} (X)$.
We start with a Young type result parallel to Lemma
\ref{inS'B}.

\par

\begin{lemma}\label{inS'}
Assume that $q,q_1,q_2\in [1,\infty ]$ satisfy \eqref{youngcond3q} or
\eqref{youngcond3q}$'$  and $s,s_1,s_2\in \mathbf R$ satisfy
\begin{equation}\label{youngcondandmore}
s_1+s_2\ge 0,\quad \text{and}\quad s\le \min (s_1,s_2 ).
\end{equation}
Then the following is true:
\begin{enumerate}
\item if \eqref{youngcond3q} holds, then the map $(f_1,f_2)\mapsto f_1
f_2$ from $\mathscr S(\rr d)\times \mathscr S(\rr d)$ to $\mathscr
S(\rr d)$ extends uniquely to a continuous mapping from $\mathscr
FL^{q_1}_{s_1}(\rr d)\times \mathscr FL^{q_2}_{s_2}(\rr d)$ to
$\mathscr FL^{q}_{s}(\rr d)$;

\vrum

\item if \eqref{youngcond3q}$'$ holds, then the map $(f_1,f_2)\mapsto
f_1 f_2$ from $C_0^\infty (X)\times C_0^\infty (X)$ to $C_0^\infty
(X)$ extends uniquely to a continuous mapping from $\mathscr
FL^{q_1}_{s_1,loc}(X)\times \mathscr FL^{q_2}_{s_2,loc}(X)$ to
$\mathscr FL^{q}_{s,loc}(X)$.
\end{enumerate}
\end{lemma}

\par

\begin{proof}
We may assume that $s= \min (s_1,s_2 )$, and
we only prove the result when $s=s_2$. The other case ($s=s_1$)
follows by similar arguments and is left for the reader. We claim that
$s_1\ge |s_2|$.

\par

In fact, this is obviously true when
$s_2\ge 0$. If instead $s_2<0$ and $s_1$ fulfills $s_1<|s_2|=-s_2$,
then $s_1+s_2<0$, which contradics \eqref{youngcondandmore}. This
proves that $s_1\ge |s_2|$, which in turn implies that if $ \omega
(\xi )=\eabs \xi ^s$, $\omega _1(\xi )= \eabs \xi^{s_1}$, and $\omega
_2(\xi )=\eabs {\xi} ^{s_2}$, then $\omega (\xi _1+\xi _2)\le C\omega
_1(\xi _1)\omega_2 (\xi _2)$ holds for some  constant $C$
which is independent of $\xi _1,\xi _2\in \rr d$.
The result now follows from Lemma \ref{inS'B}.
\end{proof}

\par

As in \cite{PTT} we denote by $L^{p,q}_1(\rr {2d})$, $p,q \in
[1,\infty]$, the mixed-norm space which consists of all $F \in L^1
_{loc} (\rr {2d} )$ such that
$$
\| F \|_{L^{p,q}_1} \equiv \Big ( \int \Big ( \int | F (\xi,
\eta)|^{p}\, d \xi\Big )^{q/p}\, d\eta \Big)^{1/q}
$$
is finite (with the usual modifications if $ p = \infty $ or if $ q =
\infty$). We also let $L^{p,q}_2(\rr {2d})$ be the set of all $F
\in L^1 _{loc} (\rr {2d} )$ such that
$$
\| F \|_{L^{p,q}_2} \equiv \Big ( \int \Big ( \int | F (\xi,
\eta)|^{q}\, d \eta \Big )^{p/q}\, d\xi \Big)^{1/p}
$$
is finite.

\par

The next result agrees with \cite[Lemma 8.3.2]{Hrm-nonlin} when $q
=2$. Here we consider extensions of the map
\begin{equation}\label{TFmap}
\begin{aligned}
(f,g)\mapsto T_F (f,g) (\xi) &= \int F(\xi,\eta) f(\eta) g (\xi-\eta) d\eta ,
\\[1ex]
&\text{from}\quad C^\infty _0 (\rr d) \times  C^\infty _0 (\rr d)  \quad \text{to}\quad  \mathscr S' (\rr d),
\end{aligned}
\end{equation}
when $F\in L^1_{loc}(\rr {2d})$.

\par

\begin{prop}\label{p-q_estimates}
Let $q \in [1,\infty]$. Then the following is true:
\begin{enumerate}
\item[(1)] if $F\in L^{\infty ,q'}_2(\rr {2d})$,
then the map in \eqref{TFmap} extends
uniquely to a continuous mapping from $ L^q (\rr d) \times L^q (\rr d)
$ to $ L^q (\rr d)$. Furthermore
\begin{equation} \label{T_F (f,g)}
\| T_F (f,g) \|_{L^q} \leq \nm {F}{L^{\infty ,q'}_2} \| f \|_{L^q} \|
g \|_{L^q}
\end{equation}
holds for every $F\in L^{\infty ,q'}_2(\rr {2d})$ and $f,g\in L^q(\rr
{d})$;

\vrum

\item[(2)] if $q>2$, $F\in L^{q,\infty}_1(\rr {2d})$,
and $r>d(1-2/q)$, then the map in \eqref{TFmap} extends uniquely to a continuous mapping from
$L^q (\rr d) \times L^q _{r}(\rr d)$ to $L^q (\rr d)$. Furthermore
\begin{equation} \label{T_F (f,g)b)}
\| T_F (f,g) \|_{L^q} \leq C\nm F{L^{q,\infty}_1} \| f \|_{L^q}
\|g \|_{L^q_r},
\end{equation}
for some  constant $C$ which is independent of $F\in
L^{q,\infty}_1(\rr {2d})$, $f\in L^q(\rr {d})$ and $g\in L^q _{r}(\rr
{d})$;

\vrum

\item[(3)] if $1\leq q \leq 2$ and $F\in
L^{q',\infty}_1(\rr {2d})$, then the map in \eqref{TFmap} extends uniquely to a continuous mapping from $ L^q (\rr
d) \times L^q (\rr d)$ to $L^q (\rr d)$. Furthermore
\begin{equation} \label{T_F (f,g)c}
\| T_F (f,g) \|_{L^q} \leq \nm {F}{L^{q',\infty }_1} \| f \|_{L^q} \|
g \|_{L^q}
\end{equation}
holds for every $F\in L^{q',\infty }_1(\rr {2d})$ and $f,g\in L^q(\rr
{d})$.
\end{enumerate}
\end{prop}

\par

\begin{proof}
The assertion (1) is obvious when $q=\infty$. Therefore assume that $q<\infty$ and that $f,g\in C_0^\infty$. By H\"older's inequality we get
\begin{multline*}
\Big (
\int | T_F (f,g)(\xi)|^q \, d\xi \Big )^{1/q}
\\[1ex]
\leq
\Big (
\int \Big ( \int | F (\xi, \eta) |^{q'} \, d\eta \Big )^{q/{q'}}
\Big (\int | f (\eta) g (\xi - \eta) |^q \, d\eta \Big )\, d\xi
\Big )^{1/q}
\\[1ex]
\leq
\nm {F}{L^{\infty ,q'}_2} \Big ( \int \int | f (\eta)|^q | g (\xi -
\eta) |^q \, d\eta d\xi \Big )^{1/q} \leq \nm {F}{L^{\infty ,q'}_2} \|
f \|_{L^q} \| g \|_{L^q}.
\end{multline*}
The result now follows from  the fact that $C_0^\infty$ is dense in $L^q$ when $q<\infty$.

\par

(2) Let $h\in C_0(\rr d)$ if $ q < \infty $ and $ h \in L^1 (\rr d)$
if $ q = \infty$.
We claim that
\begin{equation} \label{levo-desno}
\langle T_F (f,g), h \rangle =  \langle T_G (h,\check g), f \rangle,
\end{equation}
where  $G (\eta, \xi) = F (\xi, \eta) $ and $\check g (\xi ) = g
(-\xi )$.

\par

In fact, we have
\begin{multline*}
\langle T_F (f,g), h \rangle
=
\iint  F (\xi, \eta) f (\eta) g (\xi - \eta) h(\xi)\,  d\xi d\eta
\\[1ex]
= \iint  G (\xi ,\eta )  h(\eta) \check{g} (\xi - \eta ) f (\xi )\,
d\eta d\xi
\\[1ex]
= \int  T_G (h, \check{g} ) (\xi) f (\xi)\, d \xi
=   \langle T_G (h,\check g), f \rangle,
\end{multline*}
and the assertion follows.

\par

We note that $\nm {G}{L^{\infty ,q}_2}=\nm F{L^{q,\infty}_1}$. Hence
assertion (1), \eqref{T_F (f,g)} and \eqref{levo-desno} give
\begin{multline*}
|\langle T_F (f,g), h \rangle| =  | \langle T_G (h,\check g), f
\rangle | \leq \|  T_G (h,\check g) \|_{L^{q'}}  \| f \|_{L^q}
\\[1ex]
\leq  \nm {G}{L^{\infty ,q}_2} \| f \|_{L^q} \| h \|_{L^{q'}} \|  g
\|_{L^{q'}} = \nm F{L^{q,\infty}_1} \| f \|_{L^q} \| h \|_{L^{q'}} \|
g \|_{L^{q'}}.
\end{multline*}
Furthermore, by letting $q_0=q/(q-2)$, it follows that
$1/q+1/q_0=1/q'$. Hence H{\"o}lder's inequality gives
\begin{equation*}
\nm g{L^{q'}} = \nm {(g\eabs \cdo ^r)\eabs \cdo ^{-r}}{L^{q'}} \le
C\nm {g\eabs \cdo ^r}{L^q}=C\nm {g}{L^q_r},
\end{equation*}
where
$$
C= \nm {\eabs \cdo ^{-r}}{L^{q_0}}<\infty .
$$
Here the latter inequality follows from the fact that
$$
rq_0>d(1-2/q)q/(q-2)=d .
$$
The assertion (2) now follows by combining these estimates and duality.

\par

(3) Let $S$ be the trilinear map $S (F,f,g) = T_F (f,g)$, when $F\in
L^{p, \infty}_1 (\rr {2d})$, $1\le
p\le \infty$, and $f,g\in \mathscr S (\rr d)$. By straightforward
computations it follows that $S$ extends to a continuous map from
$L^{\infty, \infty}_1 (\rr {2d})\times L^1 (\rr d) \times L^1 (\rr d)$
to $ L^1 (\rr d)$ with norm $1$.

\par

Furthermore, by Cauchy-Schwartz inequality it also follows that $S$
extends to a continuous map from $L^{2, \infty}_1 (\rr {2d}) \times
L^2 (\rr d) \times L^2 (\rr d)$ to $L^2 (\rr d)$ with norm
$1$. (See also \cite[Lemma 8.3.2]{Hrm-nonlin}.)

\par

Next let $\mathcal L^{p,q}_1(\rr {2d})$ be the completion of $\mathscr
S(\rr {2d})$ under the norm $\nm \cdo {L^{p,q}_1}$. Then it follows
from Chapter 5 in \cite{BL} that complex interpolation spaces
$(\mathcal L^{p_1,q_1}_1 ,\mathcal L^{p_2,q_2}_1)_{[\theta ]}$ and
$(L^{q_1} ,L^{q_2})_{[\theta ]}$ are given by
\begin{alignat*}{3}
(\mathcal L^{p_1,q_1}_1 ,\mathcal L^{p_2,q_2}_1)_{[\theta ]} &=
\mathcal L^{p,q}_1& \quad &\text{and}&\quad
(L^{q_1} ,L^{q_2})_{[\theta ]} &= L^q
\intertext{respectively, where $p,p_j,q,q_j\in [1,\infty ]$ satisfy}
\frac{1-\theta}{p_1} +\frac{\theta}{p_2} &= \frac 1p &
\quad &\mbox{and}&\quad
\frac{1-\theta}{q_1} +\frac{\theta}{q_2} &= \frac 1q .
\end{alignat*}

\par

For each $1\le q\le 2$ it now follows by multi-linear interpolation
(cf. Theorem 4.4.1 in \cite{BL}), that
$S$ extends to a continuous map from
$\mathcal L^{q',\infty}_1 (\rr {2d}) \times L^q (\rr d)\times L^q (\rr
d)$ to $L^q (\rr d)$, with norm $1$, since this is true for $q\in \{
1,2\}$. In particular we have
\begin{equation}\label{basicnormest}
\nm {T_F(f,g)}{L^q}\le \nm {F}{L^{q',\infty} _1}\nm f{L^q}\nm g{L^q},
\end{equation}
when $F\in C_0^\infty (\rr {2d})$, $f,g\in C_0^\infty (\rr d)$ and
$1<q<2$.

\par

Now we assume that $1<q<2$ and that $F\in L^{q',\infty}_1 (\rr {2d})$
is arbitrary, and we take a sequence $F_j$ in $C^\infty _0 (\rr {2d})$
which converges to $F$ with respect to the weak$^*$ topology in
$L^{q',\infty}_1 (\rr {2d}) $ as $j$ turns to $\infty$. This is
possible since $C_0^\infty (\rr {2d}) $ is weakly dense in
$L^{q',\infty}_1 (\rr {2d}) $ when $1<q'$. Then $T_F(f,g)$ is
well-defined and smooth when $f,g\in C_0^\infty (\rr {d}) $, and since
$$
\lim _{j\to \infty} \nm {T_{F_j}(f,g)}{L^q} = \nm
{T_F(f,g)}{L^q},\quad \text{and}\quad
\lim _{j\to \infty} \nm {F_j}{L^{q',\infty} _1} = \nm F{L^{q',\infty}
_1},
$$
it follows that $T_F(f,g)\in L^q (\rr {d}) $ and that
\eqref{basicnormest} holds when $F\in L^{q',\infty}_1 (\rr {2d}) $ and
$f,g\in C_0^\infty (\rr {d}) $.

\par

The result now follows for general $F\in L^{q',\infty}_1 (\rr {2d}) $
and $f,g\in L^q (\rr {d}) $ by standard limit arguments, using the
fact that $C_0^\infty (\rr {d}) $ is dense in $L^q (\rr {d}) $ when
$1<q<2$. The proof is complete.
\end{proof}

\par

\begin{rem}
If the hypothesis in Proposition \ref{p-q_estimates} is fulfilled, and
$G(\xi ,\eta )=F(\xi ,\xi -\eta )$, then the proof of Proposition
\ref{p-q_estimates} gives that the conclusions in that proposition
still hold after $T_F$ has been replaced by $T_G$.
\end{rem}

\par

The next result agrees with Theorem 8.3.1 in \cite{Hrm-nonlin} when
$q=2$. Here we assume that $s$, $s_1$ and $s_2$ satisfy
\begin{equation}\label{s1s2scond}
0\le s_1 + s_2\quad \text{and}\quad s \le  s_1 + s_2 -
d/q',
\end{equation}
\par

\begin{thm}\label{product}
Assume that $q\in [1,\infty ]$, and let $r\ge 0$ be such that $ r >
d(1-2/q)$ when $q>2$, and $r = 0$ when  $1\le q \le 2$. Also assume
that $ s, s_j \in {\mathbf R} $ satisfy $s \leq s_j$ for $j=1,2$ and
\eqref{s1s2scond}, where the former inequality in \eqref{s1s2scond} is
strict when $s=-\min (d/q,d/q')$, and the latter inequality is strict
when $s_1=d/q'$ or $s_2=d/q'$. Then the following is true:
\begin{enumerate}
\item[(1)] the map $(f_1,f_2)\mapsto f_1 f_2$ from $\mathscr
S(\rr d)\times \mathscr S(\rr d)$ to $\mathscr S(\rr d)$ extends
uniquely to a continuous mapping from $\mathscr F L^q_{s_1} (\rr
d)\times \mathscr F L^q_{s_2+r} (\rr d)$ to $\mathscr F L^q_{s}(\rr
d)$. Furthermore,
$$
\nm {f_1f_2}{\mathscr F L^q_{s}}\le C\nm {f_1}{\mathscr F
L^q_{s_1}}\nm {f_2}{\mathscr F L^q_{s_2+r}},
$$
for some  constant $C$ which is independent of $f_1\in
\mathscr F L^q_{s_1} (\rr {d}) $ and $f_2\in \mathscr F L^q_{s_2+r}
(\rr {d}) $;

\vrum

\item[(2)] for any open set $ X$ in $\rr d$, the map $(f_1,f_2)\mapsto
f_1 f_2$ from $C_0^\infty (X)\times C_0^\infty (X)$ to $C_0^\infty
(X)$ extends uniquely to a continuous mapping from $\mathscr F
L^q_{s_1,loc} (X)\times \mathscr F L^q_{s_2+r,loc} (X)$ to
$\mathscr F L^q_{s,loc} (X)$.
\end{enumerate}
\end{thm}

\par

Since products are defined locally, Theorem \ref{product} (2) is an
immediate consequence of Theorem \ref{product} (1).

\par

We need the following lemma for the proof, concerning different
integrals of the function
\begin{equation}\label{Fdef}
F(\xi ,\eta ) = \eabs \xi ^{t_0}\eabs {\xi -\eta}^{t_1}\eabs \eta
^{t_2},
\end{equation}
where $t_0,t_1,t_2\in \mathbf R$. These integrals  are evaluated with respect to $\xi$ or $\eta$, over the sets
\begin{equation}\label{theomegasets}
\begin{aligned}
\Omega _1 &= \sets{ (\xi, \eta) \in \rr {2d}}{\eabs \eta < \delta
\eabs \xi },
\\[1ex]
\Omega _2 &= \sets{ (\xi, \eta) \in \rr {2d}}{\eabs {\xi -\eta }<
\delta \eabs \xi },
\\[1ex]
\Omega _3 &= \sets{ (\xi, \eta) \in \rr {2d}}{\delta \eabs \xi \le
\min (\eabs \eta ,\eabs {\xi -\eta }),\ |\xi |\le R},
\\[1ex]
\Omega _4 &= \sets{ (\xi, \eta) \in \rr {2d}}{\delta \eabs \xi
  \le \eabs {\xi -\eta }\le  \eabs \eta,\ |\xi |> R},
\\[1ex]
\Omega _5 &= \sets{ (\xi, \eta) \in \rr {2d}}{\delta \eabs \xi \le
  \eabs \eta \le \eabs {\xi -\eta },\ |\xi |> R},
\end{aligned}
\end{equation}
for some positive constants $\delta$ and $R$.

\par

\begin{lemma}\label{intestimates}
Assume that $F$ is given by \eqref{Fdef} and that $\Omega _1,\dots
,\Omega _5$ are given by \eqref{theomegasets}, for some
constants $0<\delta <1$ and $R\ge 4/\delta$. Also let $F_j=\chi
_{\Omega _j}F$ and $p\in [1,\infty ]$. Then, for some positive
constant $C$, the following is true:
\begin{enumerate}
\item
$$
\nm {F_1(\xi ,\cdot )}{L^p}\le
\begin{cases}
C\eabs \xi ^{t_0+t_1}\big ( 1+\eabs \xi
^{t_2+d/p} \big ),& t_2\neq -d/p ,
\\[1ex]
C\eabs \xi ^{t_0+t_1}\big ( 1+\log
\eabs \xi  \big )^{1/p},& t_2= -d/p \text;
\end{cases}
$$

\vrum

\item
$$
\nm {F_2(\xi ,\cdot )}{L^p}\le
\begin{cases}
C\eabs \xi ^{t_0+t_2}\big ( 1+\eabs \xi
^{t_1+d/p} \big ),& t_1\neq -d/p ,
\\[1ex]
C\eabs \xi ^{t_0+t_2}\big ( 1+\log
\eabs \xi  \big )^{1/p},& t_1= -d/p \text;
\end{cases}
$$

\vrum

\item $\nm {F_3(\cdo ,\eta )}{L^p}\le C\eabs \eta ^{t_1+t_2}$;

\vrum

\item if $j=4$ or $j=5$, then
$$
\nm {F_j(\cdo ,\eta )}{L^p}\le
\begin{cases}
C\eabs \eta ^{t_0+t_1+t_2+d/p},& t_0> -d/p ,
\\[1ex]
C\eabs \eta ^{t_1+t_2}\big ( 1+\log
\eabs \eta  \big )^{1/p},& t_0= -d/p,
\\[1ex]
C\eabs \eta ^{t_1+t_2},& t_0< -d/p .
\end{cases}
$$
\end{enumerate}
\end{lemma}

\par

\begin{proof}
In $\Omega _1$ we have that $\eabs \eta \le \delta \eabs \xi $,
which implies that
$$
C^{-1}\eabs \xi \le \eabs {\xi -\eta}\le C\eabs \xi
$$
for some positive constant $C$. This gives
$$
|F_1 (\xi ,\eta )|\le C\eabs \xi ^{t_0+t_1}\eabs \eta ^{t_2},
$$
for some positive constant $C$. By applying the $L^{p}$-norm with
respect to the $\eta$-variable, and using the fact that $|\eta |\le
\eabs \eta \le 2|\eta|$ when $|\eta |\ge 2$, we obtain
\begin{multline}\label{FjestA}
\nm {F_1 (\xi ,\cdot )}{L^p} \le C_1\eabs \xi ^{t_0+t_1}\Big
(\int _{|\eta|\le 2}\eabs \eta ^{t_2p}\, d\eta + \int _{2\le |\eta|
\le \eabs \xi }|\eta |^{t_2p}\, d\eta \Big )^{1/p}
\\[1ex]
\le C_2\eabs \xi ^{t_0+t_1}\Big (1 + \int _2^{\eabs \xi} r
^{t_2 p + d-1}\, dr \Big )^{1/p},
\end{multline}
for some positive constants $C_1$ and $C_2$. The result is now a
consequence of the fact that the integral on the right-hand side of
\eqref{FjestA} is estimated by
\begin{alignat*}{3}
&C(1+\eabs \xi ^{t_2+d/p}),& \quad &\text{when}& \quad t_2&\neq -d/p ,
\intertext{and}
&C(1+\log \eabs \xi )^{1/p},& \quad &\text{when}& \quad t_2&= -d/p ,
\end{alignat*}
for some positive constant $C$.

\par

The assertion (2) follows by similar arguments as in the proof of (1),
after the roles for $\eta$ and $\xi -\eta$ have been interchanged.

\par

Next we consider (3). Since $|\xi |\le R$ in $\Omega _3$, it follows
that
$$
C^{-1}\eabs \eta \le \eabs {\xi -\eta}\le C\eabs \eta ,
$$
for some positive constant $C$ which only depends on $R$. This gives
$$
F_3 (\xi ,\eta )\le C\eabs \xi ^{t_0} \eabs \eta ^{t_1+t_2},
$$
for some positive constant $C$. Since $F_3(\xi ,\eta )=0$ when $|\xi
|> R$, it follows that
$$
\nm {F_3(\cdot ,\eta)}{L^p}\le C\eabs \eta ^{t_1+t_2},
$$
for some positive constant $C$, and (3) follows.

\par

In order to estimate $\nm {F_4(\cdo ,\eta)}{L^{p}}$, we split up
$\Omega _4$ in two disjoint sets
\begin{align*}
\Omega _{4}' &= \sets {(\xi ,\eta )\in \Omega _4}{R\le |\xi | \le
|\eta |/2}
\intertext{and}
\Omega _{4}'' &= \sets {(\xi ,\eta )\in \Omega _4}{ |\eta | /2
\le |\xi |}.
\end{align*}
Then
\begin{equation}\label{F4decomp}
\nm {F_4(\cdo ,\eta)}{L^p}\le C(J_1(\eta )+J_2(\eta )),
\end{equation}
where
\begin{align*}
J_1(\eta) &= \Big (\int _{(\xi ,\eta )\in \Omega _4'}|F_4(\xi ,\eta
)|^p\, d\xi \Big )^{1/p},\quad \text{and}
\\[1ex]
J_2(\eta ) &= \Big (\int _{(\xi ,\eta )\in \Omega _4''}|F_4(\xi ,\eta
)|^p\, d\xi \Big )^{1/p}.
\end{align*}

\par

In $\Omega _4'$ we have that $|\eta |/2 \le |\xi -\eta |
\le |\eta | $. Hence, if $|\eta | \ge 2R$, then
\begin{align*}
J_1(\eta ) &\le C\Big (\int _{(\xi ,\eta )\in  \Omega _4'}\big (|\xi
|^{t_0}|\eta |^{t_1+t_2}\big )^p\, d\xi \Big )^{1/p}
\\[1ex]
&\le C\Big (\int _{R\le |\xi |\le |\eta
/2|}\big (|\xi |^{t_0}|\eta |^{t_1+t_2}\big )^p\, d\xi \Big )^{1/p}.
\end{align*}
If instead $|\eta | <2R$ we have $J_1=0$. This gives
\begin{equation}\label{J1estAA}
\begin{alignedat}{3}
J_1(\eta ) &\le C\eabs \eta^{t_0+t_1+t_2+d/p},&\quad
&\text{when}&\quad t_0 &>-d/p,
\\[2ex]
J_1(\eta ) &\le C\eabs \eta ^{t_1+t_2}\big (\log |\eta
|)^{1/p},&\quad &\text{when}&\quad t_0 &=-d/p
\\[2ex]
J_1(\eta ) &\le C\eabs \eta ^{t_1+t_2},&\quad &\text{when}&\quad
t_0 &<-d/p.
\end{alignedat}
\end{equation}

\par

For $\Omega _4''$ we claim that $|\xi |\le R|\xi -\eta |$. Admitting
this for a while it follows from the assumptions that
$$
\frac 12\le \frac {|\eta |}2\le |\xi | \le R|\xi -\eta | \le R|\eta
|.
$$
By the first part of the proof of (4), it follows that
\eqref{J1estAA} holds after $J_1(\eta )$ is replaced by $J_2(\eta
)$. The assertion (4) in the case $j=4$ now follows by combining these
estimates with \eqref{F4decomp}.

\par

It remains to prove $|\xi |\le R|\xi -\eta |$. By the assumptions we
have
\begin{equation}\label{moreest22}
\frac {16}{\delta ^2}\le |\xi |^2\le \frac {4}{\delta ^2}-1+\frac
{4}{\delta ^2}|\xi -\eta |^2,
\end{equation}
which implies that $|\xi -\eta|\ge 1$. Hence, \eqref{moreest22} gives
$$
|\xi |^2\le \frac {4}{\delta ^2}-1+\frac {4}{\delta ^2} |\xi -\eta |^2 \le \frac {16}{\delta ^2}|\xi -\eta |^2 =R^2|\xi -\eta |^2.
$$
This proves that $|\xi |\le R|\xi -\eta |$.

\par

The assertion (4) in the case $j=5$ follows by similar arguments as
for $j=4$, after the roles of $\eta$ and $\xi -\eta$ have been
interchanged. The proof is complete.
\end{proof}

\par

\begin{proof}[Proof of Theorem \ref{product}.]
Let $t_0=s$, $t_1=-s_1$, $t_2=-s_2$, and let $F$, $F_j$ and $\Omega
_j$ be the same as in Lemma \ref{intestimates} for $j=1,\dots ,5$,
after $\Omega _2$ has been modified into
$$
\Omega _2 = \sets{ (\xi, \eta) \in \rr {2d}}{\eabs {\xi -\eta }<
\delta \eabs \xi }\setminus \Omega _1,
$$
and let $T_F$ be the same as in Proposition \ref{p-q_estimates}.
Then $\cup \Omega _j=\rr {2d}$, $\Omega _j\cap
\Omega _k$ has Lebesgue measure zero when  $j\neq k$, and
\begin{multline*}
\eabs \xi ^s \mathscr F(f_1f_2) (\xi) = (2\pi)^{-d/2}\eabs \xi ^s
\widehat f_1 * \widehat f_2 (\xi)
\\[1ex]
= (2\pi)^{-d/2} \int F (\xi, \eta ) u_1 (\xi - \eta) u_2 (\eta) d\eta
=T_F(u_1,u_2)
\\[1ex]
=T_{F_1}(u_1,u_2)+\cdots + T_{F_5}(u_1,u_2)
\end{multline*}
where $ u_j (\xi) = \eabs \xi ^{s_j} \widehat f_j  (\xi)$, $ j=1,2$.

\par

By Lemma \ref{intestimates}
(1) and (2) it follows that
\begin{equation}\label{Fjest}
\nm {F_j}{L^{\infty ,q'}_2}<\infty ,\qquad j=1,2
\end{equation}
when $j=1,2$. Furthermore, by  Lemma \ref{intestimates} (3)--(4) we
get
\begin{equation}\label{Fjest2}
\nm {F_j}{L^{\infty ,p}_1}<\infty ,\qquad j=3,4,5,
\end{equation}
when $p=\max (q,q')$. Hence, Proposition \ref{p-q_estimates} (1),
applied to $T_{F_j}(u_1,u_2)$ when $ j = 1, 2 $ and  Proposition
\ref{p-q_estimates} (2) and (3), applied to
$T_{F_j}(u_1,u_2)$ when $ j = 3, 4, 5 $,
show that
$T_{F_j}(u_1,u_2)$ is well-defined when $u_1\in L^q (\rr {d}) $ and
$u_2\in L^q_{r} (\rr {d}) $, and that
\begin{equation}\label{TFjest}
\nm {T_{F_j}(u_1,u_2)}{L^q}\le C \nm {u_1}{L^q}\nm {u_2}{L^q_r}\qquad
j=1,\dots ,5.
\end{equation}

\par

The result now follows from \eqref{TFjest} and the following
relations
\begin{align*}
\nm {T_{F}(u_1,u_2)}{L^q} &= \nm {f_1f_2}{\mathscr FL^q_s},\quad
\nm {u_1}{L^q} = \nm {f_1}{\mathscr FL^q_{s_1}}
\\[1ex]
\text{and}\quad  \nm {u_2}{L^q_r} &= \nm {f_1}{\mathscr FL^q_{s_2+r}}.
\end{align*}
The proof is complete.
\end{proof}

\par

We have now the following.

\par

\begin{prop}\label{algebramodules}
Assume that $q,q_0\in [1,\infty ]$ are such that $q_0\le q$, $s\ge
d/q'$ when $1\le q<2$ and $s>d(3/q'-1)$ when $2\le q\le \infty$, and that
$X\subseteq \rr d$ is open. Then the following is true:
\begin{enumerate}
\item $(\mathscr FL^q_s(\rr d),\cdo )$ is an algebra, and a
(left-right) $\mathscr FL^{q_0}_s(\rr d)$-module. Furthermore, if
$f_1,\dots ,f_N\in \mathscr FL^q_s(\rr d)$ and $g\in \mathscr
FL^{q_0}_s(\rr d)$, then
$$
\nm {f_1\cdots f_N\cdot g}{\mathscr FL^q_s}\le C^{N+1}\nm
{f_1}{\mathscr FL^q_s}\cdots \nm {f_N}{\mathscr FL^q_s}\cdot \nm
{g}{\mathscr FL^{q_0}_s}\text ;
$$

\vrum

\item $(\mathscr FL^q_{s,loc}(X),\cdo )$ is an algebra, and a
(left-right) $\mathscr FL^{q_0}_{s,loc}(X)$-module.
\end{enumerate}
\end{prop}

\par

\begin{proof}
The assertion (1) in the case $q_0\in \{ 1,q\}$, follows immediately
from Lemma \ref{inS'} and Theorem \ref{product}. The result now
follows for general $q_0$ by interpolation.

\par

Since the product is a local operation, (2) is an immediate
consequence of (1). The proof is complete.
\end{proof}

\par

\section{Micro-local properties in the Fourer Lebesgue spaces}
\label{sec4}

\par

In this section we discuss wave-front properties of Fourier Lebesgue
types for products of distributions.

\par

We start with the micro-local characterization of the product
$f_1f_2$ when $f_j \in \mathscr F L^q_{s_j,loc}(X)$, $j=1,2$.

\par

\begin{thm}\label{microlocal-functions}
Assume that $q \in [1,\infty]$ and let $f_j \in  \mathscr F
L^q_{s_j,loc} (X)$, $j = 1,2$. Then the following is true:
\begin{enumerate}
\item if $s_1 - |s_2| \ge 0$ when $q=1$ and $  s_1 - |s_2| > d/q'$ otherwise, then
$$
WF_{\mathscr F L^{q} _{s_2}}(f_1f_2) \subseteq WF_{\mathscr F
L^q_{|s_2|}} (f_1)\text ;
$$

\vrum

\item if instead $s_1 + s_2 \geq s \geq 0$ when $q=1$ and $ s_1 + s_2 - d/q' > s \geq 0$ otherwise, and
$ s_2 - s \ge d/q' $, then
$$
WF_{\mathscr F L^{q} _{s}}(f_1f_2) \subseteq WF_{\mathscr F
L^q_{|s_2|}} (f_1)\text .
$$
\end{enumerate}
\end{thm}

\par

We note that $f_1 f_2$ in Theorem \ref{microlocal-functions} makes
sense as an element in $\mathscr FL^q_{s,loc} (X)$, for some $s\in
\mathbf R$.

\par

For the proof we need the following result parallel to Lemma
\ref{intestimates}.

\par

\begin{lemma}\label{intestimates2}
Assume that $F$ is given by \eqref{Fdef} and
$$
\Omega = \sets {(\xi ,\eta)\in \rr {2d}}{|\eta -\xi|\ge c|\xi |,\
\eabs \xi \ge R,\ \eabs \eta \ge R},
$$
for some $c>0$ and $R>1$. Also let $p\in [1,\infty ]$ be such that
$t_1+t_2<-d/p$ when $p<\infty$ and $t_1+t_2\le 0$ when $p=\infty$.
Then
$$
\nm {F(\xi ,\cdot )}{L^p(\Omega )}\le
\begin{cases}
C\eabs \xi ^{t_0},& t_2\geq -d/p ,
\\[1ex]
C\eabs \xi ^{t_0}(1+\eabs \xi ^{t_1}),& t_2< -d/p ,
\end{cases}
$$
for some constant $C>0$.
\end{lemma}

\par

\begin{proof}
Again we only prove the result in the case $p<\infty$. The simple
modifications to the case $p=\infty$ are left for the reader. Let
\begin{align*}
\Omega _1 &= \sets {(\xi ,\eta)\in \Omega}{|\eta -\xi|\le C|\xi |}
\\[1ex]
\Omega _2 &= \sets {(\xi ,\eta)\in \Omega}{|\eta -\xi|\ge C|\xi |}
\end{align*}
for some constant $C>3$. From the assumptions it follows that
\begin{equation}\label{xietacond}
\begin{aligned}
c\le |\xi |/\eabs \xi \le 1,\quad c&\le |\xi -\eta |/\eabs {\xi
-\eta} &\le 1,
\\[1ex]
\text{and}\quad c &\le |\eta |/\eabs \eta \le 1,
\end{aligned}
\end{equation}
for some constant $c>0$. This gives
\begin{equation}\label{FOmegaest2}
\nm {F(\xi ,\cdot )}{L^p(\Omega )}\le C\eabs \xi ^{t_0}(I_1+I_2)
\end{equation}
for some positive constant $C$, where
\begin{equation*}
I_j^p = \int _{(\xi ,\eta )\in \Omega _j}|\eta - \xi |^{t_1p}|\eta
|^{t_2p}\, d\eta ,\quad j=1,2.
\end{equation*}

\par

We need to estimate $I_1$ and $I_2$, and start by considering
$I_1$. Since $c|\xi |\le |\eta -\xi |\le C|\xi |$ and $|\eta |\ge
C$ in $\Omega _1$, it follows that $C\le |\eta |\le C_1|\xi |$ for
some constant $C_1$, and
\begin{multline*}
I_1^p \le C_2|\xi |^{t_1p} \int _{C\le |\eta |\le C_1|\xi |}|\eta
|^{t_2p}\, d\eta
\\[1ex]
= C_3|\xi |^{t_1p} \int _C^{C_1|\xi |}r^{t_2p+d-1}\, dr \le
\begin{cases}
C_4 |\xi |^{(t_1+t_2)p+d},& t_2> -d/p ,
\\[1ex]
C_4 |\xi |^{t_1}(1+\log \eabs \xi ),& t_2 = -d/p ,
\\[1ex]
C_4 |\xi |^{t_1p},& t_2 < -d/p,
\end{cases}
\end{multline*}
for some positive constants $C,C_1,\dots C_4$. A combination of
\eqref{xietacond} and the facts $t_1<0$ when $t_2=-d/p$ and
$(t_1+t_2)p+d<0$ now give
\begin{equation}\label{I1est1}
I_1 \le
\begin{cases}
C ,& t_2 \ge -d/p,
\\[1ex]
C \eabs {\xi}^{t_1},& t_2 < -d/p,
\end{cases}
\end{equation}
for some positive constant $C$.

\par

It remains to estimate $I_2$. From the assumptions it follows that
$2^{-1}|\eta |\le |\eta -\xi|\le 2|\eta |$ when $(\xi ,\eta )\in
\Omega _2$. This implies that
\begin{equation}\label{I2est1}
I_2^p \le C\int _{|\eta |\ge c}|\eta |^{(t_1+t_2)p}\, d\eta
<\infty .
\end{equation}
The result now follows by combining \eqref{FOmegaest2},
\eqref{I1est1} and \eqref{I2est1}.
\end{proof}

\par

\begin{proof}[Proof of Theorem \ref{microlocal-functions}]
Again we only prove the result for $1<q<\infty$. The small
modifications to the cases $q=1$ or $q=\infty$ are left for the
reader. We shall mainly follow the ideas in the proof of
\cite[Theorem 8.3.3]{Hrm-nonlin}.

\par

Since the statements only concern local properties, we may assume that
$f_1\in \mathscr FL^q_{s_1} (\rr d)$ and $f_2\in \mathscr FL^q_{s_2}
(\rr d)$ are compactly supported.

\par

(1) Assume that $(x_0,\xi _0)\notin WF_{\mathscr FL^q _{|s_2|}}
(f_1)$. Then $|f_1|_{\mathscr FL^{q,\Gamma}_{|s_2|}}<\infty$ for
some $\Gamma =\Gamma _{\xi _0}$. It is sufficient to prove that
$|f_1f_2|_{\mathscr FL^{q,\Gamma _1}_{s_2}} <\infty$ when $\overline
{\Gamma _1}\subseteq \Gamma$.

\par

Since $\eabs \xi ^{s_2}\le \eabs {\xi -\eta}^{s_2}\eabs \eta
^{|s_2|}$ we get
$$
\eabs \xi ^{s_2}|\mathscr F(f_1f_2)(\xi )|\le C(J_1(\xi )+J_2(\xi ))
$$
for some positive constant $C$, where
\begin{alignat*}{2}
J_1(\xi ) &= \eabs \xi ^{s_2} (|\widehat f_2|*|\chi _{\complement
\Gamma}& & \widehat f_1|)(\xi )
= \eabs \xi ^{s_2}\int_{\complement \Gamma}
|\widehat f_2 (\xi - \eta)| |\widehat f_1 (\eta)|\, d\eta
\intertext{and}
J_2(\xi ) &= (|\widehat f_2\eabs \cdo ^{s_2}|*|\chi _{\Gamma}& &\widehat
f_1\eabs \cdo ^{|s_2|} |)(\xi )
\\[1ex]
& & & = \int_{\Gamma} |\widehat f_2 (\xi
-\eta)\eabs{\xi -\eta}^{s_2}|\,  |\widehat f_1 (\eta)\eabs
\eta ^{|s_2|}|\, d\eta .
\end{alignat*}
We need to estimate $\nm {J_1}{L^q(\Gamma _1)}$ and $\nm
{J_2}{L^q(\Gamma _1)}$.

\par

By Minkowski's inequality, we get
\begin{multline} \label{zvezdica}
\nm {J_2}{L^q(\Gamma _1)} = \Big ( \int _{\Gamma _1}\Big ( \int
_{\Gamma} |\widehat f_2(\xi -\eta )|\eabs {\xi -\eta} ^{s_2} |\widehat
f_1(\eta )|\eabs \eta ^{|s_2|}\, d\eta \Big )^q\, d\xi \Big )^{1/q}
\\[1ex]
\le \int _{\Gamma} \Big ( \int _{\Gamma _1} \Big ( |\widehat f_2(\xi
-\eta )|\eabs {\xi -\eta} ^{s_2} |\widehat f_1(\eta )|\eabs \eta
^{|s_2|} \Big )^q \, d\xi \Big )^{1/q} \, d\eta
\\[1ex]
\le \nm {\widehat f_1\eabs \cdot ^{|s_2|}}{L^1(\Gamma )}
\nm { \widehat f_2 \langle \cdot \rangle ^{s_2}}{L^q(\rr d)} =
|f_1|_{\mathscr F L^{1,\Gamma }_{|s_2|}}\nm {f_2}{\mathscr
FL_{s_2}^q}.
\end{multline}

\par

Now we have
$$
|f_1|_{\mathscr F L^{1,\Gamma }_{|s_2|}}\le
\nm {\widehat f_1 \eabs \cdot ^{s_1} \eabs \cdot ^{|s_2| - s_1}}{L^1}
\leq
\nm {\widehat f_1 }{L^q _{s_1}}
\nm { \eabs \cdot ^{|s_2| - s_1}}{L^{q'}}
\leq C \| f_1 \|_{L^q _{s_1}},
$$
since $ (s_1 - |s_2|) q' > d$. Hence $\nm
{J_2}{L^q(\Gamma_1)}<\infty$.

\par

The assertion (1) therefore follows if we prove $\nm {J_1}{L^q(\Gamma
_1)}<\infty$. If $ \xi \in \Gamma _1$
and $ \eta \not\in \Gamma$, then
it follows that $ \eabs \xi \leq C \langle \xi - \eta \rangle$ for
some  constant $C$. Hence, if $s_2 \geq 0$, we obtain
$$
|J_1 (\xi )| \leq C\int_{\complement \Gamma}
|\widehat f_1 (\eta)| \, |\widehat f_2 (\xi - \eta)|\eabs{\xi - \eta}
^{s_2}\, d\eta, \quad \xi \in \Gamma _1,
$$
for some  constant $C$. Since $f_1 \in \mathscr FL^1 $ because
$s_1 > d/q'$, there exists a constant $C>0$ such that
$$
\nm {J_1}{L^q(\Gamma _1)} \leq  C \| f_1 \|_{\mathscr FL^1} \|
f_2\|_{\mathscr FL^q_{s_2}} < \infty ,
$$
and (1) follows in the case $s_2\ge 0$.

\par

Assume instead that $s_2 < 0$. Then,
\begin{multline*}
|J_1 (\xi)| \leq \int_{\complement \Gamma}
|\widehat f_1 (\eta)|\eabs \eta ^{-s_2}   |\widehat f_2 (\xi -
\eta)|\eabs{\xi - \eta} ^{s_2}\,  d\eta
\\[1ex]
\leq
\big (\, ( | \widehat f_1 |\eabs \cdot ^{-s_2}) *(|\widehat f_2 |\eabs
\cdot ^{s_2} )\, \big ) (\xi),
\end{multline*}
and Minkowski's inequality gives
$$
\nm {J_1}{L^q(\Gamma _1)} \leq \nm {f_1}{\mathscr FL^1_{-s_2}} \nm
{f_2}{\mathscr FL^q_{s_2}}.
$$
By H{\"o}lder's inequality we get
$$
\nm {f_1}{\mathscr FL^1_{-s_2}} \le C\nm {f_1}{\mathscr FL^q_{s_1}},
$$
where $C=\nm {\eabs \cdot ^{-(s_1+s_2)}}{L^{q'}}$ is finite since $
s_1 + s_2 > d/q'$. Summing up we have proved
$$
\nm {J_1}{L^q(\Gamma _1)} \leq C\nm {f_1}{\mathscr FL^q_{s_1}}\nm
{f_2}{\mathscr FL^q_{s_2}}<\infty .
$$
The estimates of $ \nm {J_2}{L^q(\Gamma _1)}$ are proved in a similar way as for $\nm {J_1}{L^q(\Gamma _1)}$. This proves (1).

\par

(2) Assume that $(x_0,\xi _0)\notin WF_{\mathscr FL^q _{|s_2|}} (f_1)$,
and let $J_1$, $J_2$, $ \Gamma $ and $\Gamma_1 $ be the same as in the
first part of the proof, after $s_2$ has been replaced by $s$. By the
assumptions we have $0\le s<s_2$. Hence, by (\ref{zvezdica}) with $s$
instead of $s_2$ it follows that
\begin{multline*}
\nm {J_2}{L^q(\Gamma _1)} \le |f_1| _{L^{q }_{s} (\Gamma _1)}\nm
{f_2}{\mathscr FL_{s}^1} \le |f_1| _{L^{q}_{s} (\Gamma _1)} \nm
{f_2}{\mathscr F L^1 _{s}}
\\[1ex]
\leq C |f_1| _{L^{q}_{s} (\Gamma _1)} \nm {f_2}{\mathscr F L^q _{s_2}}
<\infty .
\end{multline*}

\par

It remains to prove $\nm {J_1}{L^q(\Gamma _1)}<\infty$. For $ \xi \in
\Gamma _1$ and $ \eta
\not\in \Gamma$, we have $ \langle \xi \rangle ^s \leq C \langle  \xi -
\eta \rangle ^s $, since $ s\geq 0 $. This gives
\begin{equation*}
J_1(\xi )
\leq
\int_{\rr d} F (\xi, \eta)|\widehat f_2 (\xi - \eta)\langle \xi - \eta
\rangle ^{s_2} |\,
|\widehat f_1 (\eta) \eabs \eta ^{s_1}|\,  d\eta
\end{equation*}
where
$$
F (\xi, \eta ) =
\left \{
\begin{array}{lll}
C \eabs \eta ^{-s_1}  \langle \xi - \eta \rangle ^{s - s_2}, &
\mbox{when} & \eta \in \complement \Gamma ,\ \ \xi \in \Gamma _1
\\[1ex]
0, & \mbox{otherwise.} &
\end{array}
\right .
$$

\par

The assertion follows from Proposition \ref{p-q_estimates} (1) if we show
that $h_\Omega \in L^\infty$ with $\Omega =\rr d$, where $h_\Omega (\xi) \equiv \nm {F(\xi ,\cdo )}{L^{q'}(\Omega )}$.
It is obvious that $h_{\Omega}(\xi )$ with $\Omega =\rr d$ is bounded when $\xi$ belongs to a
bounded set, since $s_1 + s_2 - s > d/q'$. Furthermore, by the
assumptions it follows that $s_2>s$, which implies that if $K\subseteq
\rr d$ then $h_K\in L^\infty $. Therefore, the result
follows if we prove that $h_\Omega (\xi )$ is bounded when $|\xi |\ge C$ and $\Omega ={\sets {\eta \in \rr d}{|\eta |\ge C}}$, for some large positive constant $C$.
$$
h_0(\xi) \equiv \nm {F_0(\xi ,\cdo )}{L^{q'}}.
$$
Here $F_0(\xi ,\eta )=F(\xi ,\eta )$ when $|\xi |\ge C$ and $|\eta
|\ge C$ for some large constant $C$, and $F_0(\xi ,\eta )=0$
otherwise. This follows by taking $p=q'$, $t_0=0$, $t_1=s-s_2$ and
$t_2=-s_1$ in Lemma \ref{intestimates2}. The proof is complete.
\end{proof}

\par

In the next theorem we consider the "critical case"
$s = s_1 + s_2 - \min (d/q,d/q')$ comparing to Theorem \ref{microlocal-functions} (2).

\par

\begin{thm}\label{microlocal-functions-some-more}
Assume that $q \in [1,\infty]$, $r=0$ when $1\le q\le 2$, $r>d(1-2/q)$
when $q>2$, and that $s,s_j,N_j\in \mathbf R$, $j=1,2$, satisfy
\begin{gather}
s_1 + s_2 > 0,\qquad s = s_1+s_2-\min (d/q,d/q'),\notag
\\[1ex]
\begin{aligned}\label{N1N2cond}
N_1 &\ge s_1+|s_2| + \max (0,d(1-2/q)),\quad \text{and}
\\[1ex]
N_2 &\ge s_2+|s_1| + \max (0,d(1-2/q)),
\end{aligned}
\end{gather}
with strict inequalities in \eqref{N1N2cond} when $q<\infty$. If $f_1
\in  \mathscr F L^q_{s_1,loc} (X)$ and $f_2 \in \mathscr F
L^q_{s_2+r,loc} (X)$, then
$$
WF_{\mathscr F L^{q} _{s}}(f_1 f_2) \subseteq WF_{\mathscr F
L^q_{N_1}} (f_1)\cup WF_{\mathscr F L^q_{N_2}} (f_2)\text .
$$
\end{thm}

\par

\begin{proof}
Again we only prove the result for $1<q<\infty$, leaving the small
modifications when $q\in \{ 1,\infty \}$ to the reader.

\par

Assume that $(x_0,\xi _0)\notin WF_{\mathscr FL^q_{N_j}}(f_j)$, $j
= 1,2$. It is no restriction to assume that $f_j$ has compact support
and $\xi _0\notin \Sigma_{\mathscr FL^q_{N_j}}(f_j)$. Then
$|f_j|_{\mathscr FL^{q,\Gamma}_{N_j}}<\infty$ for some cone $\Gamma$
of $\xi _0$. Furthermore, for some $\delta \in (0,1)$ and open cone $
\Gamma _1$ of $\xi _0$ such that $\overline  \Gamma _1\subseteq \Gamma
$ we have $\xi - \eta \in \Gamma$ when $\xi \in \Gamma _1$ and $|\eta
| < \delta |\xi |$. Let $\Omega _1$ and $\Omega _2$ be the same as in
\eqref{theomegasets}, $\Omega _0=\complement (\Omega _1\cup \Omega
_2)$, and let
$$
J_k(\xi ) = \int _{(\xi ,\eta )\in \Omega _k}\eabs \xi ^s
|\widehat f_2(\xi -\eta )\widehat f_1(\eta )|\, d\eta
$$
for $k=0,1,2$. The result follows if we prove that
\begin{align}
\nm {J_0}{L^q(\Gamma _1)} &\le C\nm {f_1}{\mathscr FL^q_{s_1}}\nm
{f_2}{\mathscr FL^q_{s_2+r}}\label{J0est},
\\[1ex]
\nm {J_1}{L^q(\Gamma _1)} &\le C\nm {f_1}{\mathscr FL^q_{s_1}}
|f_2|_{\mathscr FL^{q,\Gamma}_{N_2}},\label{J1est}
\intertext{and}
\nm {J_2}{L^q(\Gamma _1)} &\le C|f_1|_{\mathscr FL^{q,\Gamma}_{N_1}}
\nm {f_2}{\mathscr FL^q_{s_2}}
,\label{J2est}
\end{align}
for some constant $C>0$.

\par

In order to prove \eqref{J1est}, we choose $\ep >0$,
$N_0$ and $N$ such that
$$
N_0=s_1+s_2+\ep,\quad N=N_0-s_1+|s_1|\quad \text{and}
\quad N < N_2 - d|1-2/q|.
$$

\par

We have $C^{-1}\eabs \xi \le \eabs {\xi -\eta}\le C\eabs \xi$ when
$(\xi ,\eta )\in \Omega _1$. This gives
\begin{multline} \label{J1estC}
J_1(\xi ) = \eabs \xi ^{s-{N_0}} \int_{\Omega _1}
\eabs \xi  ^{N_0} |\widehat f_2 (\xi - \eta)|\, |\widehat f_1 (\eta)
\eabs \eta  ^{s_1} | \eabs \eta ^{-s_1}\, d\eta
\\[1ex]
\leq
\eabs \xi ^{s-N_0} \int_{\Omega _1} \eabs \xi ^{N_0 -s_1}
|\widehat f_2 (\xi - \eta) \langle \xi -  \eta \rangle ^{|s_1|}|\, |
\widehat f_1 (\eta) \langle  \eta \rangle ^{s_1}|\, d\eta
\\[1ex]
\leq
C \langle \xi \rangle^{s-N_0} \int_{\Omega _1}  |\widehat f_2 (\xi -
\eta)\langle \xi -  \eta \rangle ^{N_0 -s_1 +|s_1|}|\, |\widehat f_1
(\eta) \langle  \eta \rangle ^{s_1}|\, d\eta
\end{multline}
A combination of \eqref{J1estC} and H{\"o}lder's inequality give
$$
J_1(\xi )\le \langle \xi \rangle^{-\ep -d/q} \nm {f_1}{\mathscr
FL^{q}_{s_1}}|f_2|_{\mathscr FL^{q',\Gamma }_{N}}.
$$
By applying the $L^q(\Gamma _1)$ norm we obtain
\begin{equation}\label{J1estno2}
\nm {J_1}{L^q(\Gamma _1)} \le C\nm {f_1}{\mathscr FL^q_{s_1}}|f_2
|_{\mathscr FL^{q',\Gamma }_{N}}.
\end{equation}
Hence \eqref{J1est} follows from \eqref{J1estno2} if we prove that
\begin{equation}\label{estf2againA}
|f_2 |_{\mathscr FL^{q',\Gamma }_{N}}\le C|f_2|_{\mathscr
FL^{q,\Gamma}_{N_2}},
\end{equation}
for some  constant $C>0$.

\par

If $1\le q\le 2$, then \eqref{estf2againA} follows from the facts that
$f_2$ has compact support, $N\le N_2$ and $q\le q'$. Therefore, assume
that $q>2$, and set $q_0=q/(q-2)$. Then $q_0\ge 1$ and
$1/q'=1/q+1/q_0$. Let $\ep _0>0$. By H{\"o}lder's
inequality we obtain
\begin{multline}\label{f2estimatesA}
|f_2 |_{\mathscr FL^{q',\Gamma }_{N}} = \nm {\widehat f_2\cdot \eabs
\cdo ^{N}}{L^{q'}(\Gamma )}
\\[1ex]
=  \nm {(\widehat f_2\cdot \eabs \cdo ^{N+(d+\ep _0)/q_0})\eabs \cdo
^{-(d+\ep _0)/q_0}}{L^{q'}(\Gamma )}
\\[1ex]
\le C\nm {\widehat f_2\cdot \eabs \cdo ^{N+(d+\ep _0)/q_0}}{L^{q}(\Gamma
)} =C|f_2 |_{\mathscr FL^{q,\Gamma }_{N+(d+\ep _0)/q_0}},
\end{multline}
where $C=\nm {\eabs \cdo ^{-d-\ep_0}}{L^1}^{1/q_0}<\infty$. By choosing
$\ep _0$ small enough, it follows that $N+(d+\ep _0)/q_0$ is smaller than
$N_2$. Hence \eqref{estf2againA} follows from
\eqref{f2estimatesA}. This proves \eqref{J1est}, and
the estimate \eqref{J2est} follows by similar arguments, after the
roles of $f_1 $ and $ f_2 $ have been interchanged. The details are
left for the reader.

\par

It remains to prove \eqref{J0est}. Let $t_0=s$,
$t_1=-s_1$, $t_2=-s_2$, and let $\Omega _j$ and $F_j$ for
$j=3,4,5$ be the same as in Lemma \ref{intestimates}. Also let $T_F$
be the same as in Proposition \ref{p-q_estimates}. Since
$$
\sets {(\xi ,\eta )\in \Omega _0}{\xi \in \Gamma _1}\subseteq \Omega
_3\cup \Omega _4\cup \Omega _5,
$$
it follows that
$$
J_0\le T_{F_3}(u_1,u_2)+T_{F_4}(u_1,u_2)+T_{F_5}(u_1,u_2),
$$
where $u_j(\xi )=|\widehat f_j(\xi )\eabs \xi ^{s_j}|$. Hence it
suffices to prove that
$$
\nm {T_{F_j}(u_1,u_2)}{L^q(\Gamma _1)}\le C\nm {f_1}{\mathscr
FL^q_{s_1}}\nm {f_2}{\mathscr FL^q_{s_2+r}}.
$$

\par

By Lemma \ref{intestimates} and the
assumptions on $s$, $s_1$ and $s_2$, it follows that $\nm
{F_j}{L^{q,\infty}_1}<\infty$ when $j=3,4,5$.
Hence Proposition \ref{p-q_estimates} (2) gives
$$
\nm {T_{F_j}(u_1,u_2)}{L^q(\Gamma _1)}\le C\nm {u_1}{L^q}\nm
{u_2}{L^q_r} = C\nm {f_1}{\mathscr FL^q_{s_1}}\nm {f_2}{\mathscr
FL^q_{s_2+r}}
$$
when $j = 3,4,5$, and \eqref{J0est} follows. The proof is complete.
\end{proof}

\par

\section{Semi-linear equations} \label{sec5}

\par

In this section we apply the wave-front results of previous sections
to solutions of a broad class of semi-linear differential equations,
and more general algebraic expressions of distributions.

\par

These expressions are of the form
\begin{equation}\label{gpoldef}
\mathcal G(x)\equiv G(x,f_1(x),\dots ,f_N(x))
\end{equation}
when $f_1,\dots ,f_N$ are appropriate
distributions and $G(x,y)$ is an appropriate polynomial in the
$y$-variable, i.{\,}e. $G$ is of the form
\begin{equation}\label{velikoF}
G(x,y) = \sum _{0<|\alpha |\le m}a_\alpha (x)y^\alpha ,
\end{equation}
where $a_\alpha$ are appropriate distributions and  $ y = (y_1, \dots,
y_N) $.

\par

In the following result we assume that $a_\alpha$, $ 0<|\alpha |\le m$,  locally belong to
appropriate classes of Fourier Lebesgue spaces.

\par

\begin{thm}\label{F}
Assume that $X\subseteq \rr d$ is open, $q\in [1,\infty]$ and $ r
\geq d/q' $ with the strict inequality when $q=\infty$. Let
$$
s \geq d/q',\quad s\leq \sigma \leq 2s - d/q',
$$
and let $G$ and $\mathcal G$ be given by \eqref{gpoldef} and \eqref{velikoF}
for some integer $m> 0$, where $a_\alpha \in
\mathscr FL^1_{\sigma ,loc}(X)$, $ 0<|\alpha |\le m$,  and $f_j \in \mathscr
FL^q_{s+(m-1)r,loc}(X)$, $ j = 1,\dots,N$. Then the following is
true:
\begin{enumerate}
\item[(1)]  $\mathcal G$ in \eqref{gpoldef} makes sense as an element in
$\mathscr FL^q_{s,loc}(X)$;

\vrum

\item[(2)] $  WF_{\mathscr F L^{q} _{\sigma}}(\mathcal G) \subset
\cup_{j =1} ^{N} WF_{\mathscr FL^{q} _{\sigma +(m-1)r}} (f_j)$.
\end{enumerate}
\end{thm}

\par

\begin{proof}
Assertion (2), in the case $m=1$, and assertion (1) follow
immediately from Proposition \ref{algebramodules}.

\par

We need to prove (2) for $m\ge 2$. By Proposition
\ref{algebramodules} we may assume that $a_\alpha$ are constants, $ 0<|\alpha |\le m$.

\par

Assume that $m=2$ and that and $ (x_0 , \xi_0 ) \not\in
WF_{\mathscr FL^{q} _{\sigma + r}} (f_j)$, $j=1,\dots ,N$. Then
each term in \eqref{gpoldef} is of the form
\begin{equation}\label{Gterms}
a_jf_j\quad \text{or} \quad a_{j,k}f_jf_k,\quad j,k=1,\dots ,N
\end{equation}
We write $ f_j= g_j + h_j$, $j=1,\dots ,N$ where $ h_j \in
\mathscr F L^{q} _{\sigma + r, loc} (X)$ and  $ (x_0 , \xi_0 )
\not\in WF (g_j)$, in view of Proposition \ref{opadanje-za-sve-s}.
Then $ g_j \in \mathscr F L^{q} _{s+ r, loc}(X)$ by the
assumptions.

\par

First we consider the case $q\geq 2$. We need to prove that
\begin{equation}\label{WFfjfk}
(x_0 , \xi_0 ) \not\in WF_{ \mathscr F L^{q} _{\sigma}} (f_j f_k),
\end{equation}
and then it suffices to prove that
\begin{align}
(x_0 , \xi_0 ) \not\in WF_{ \mathscr F L^{q} _{\sigma}} (g_j
g_k),\quad (x_0 , \xi_0 ) &\not\in WF_{ \mathscr F L^{q}
_{\sigma}} (g_j h_k),\label{decompfjfk1}
\\[1ex]
\text{and}\qquad (x_0 , \xi_0 ) &\not\in WF_{ \mathscr F L^{q}
_{\sigma}} (h_j h_k).\label{decompfjfk2}
\end{align}

\par

For \eqref{decompfjfk2} we note that $ r
d(1 - 2/q) $ since $ d/q' > d(1 - 2/q)$. Hence, by Theorem
\ref{product} (2) we have $ h_j h_k \in \mathscr F L^{q} _{\sigma
, loc} (X)$, and \eqref{decompfjfk2} follows.

\par

In order to prove \eqref{decompfjfk1} we set $ s_1 = \sigma + r $,
$ s_2 = \sigma$ and $ s = \sigma$. Then $s_1 - s = r \geq d/q'$,
and by Theorem \ref{microlocal-functions} (2) it follows that $
WF_{ \mathscr F L^{q} _{\sigma } } (g_j g_k) \subseteq WF_{
\mathscr F L^{q} _{\sigma + r} } (g_j) $. This gives the first
relation in \eqref{decompfjfk1}, and the second one follows by the
same theorem. This proves \eqref{WFfjfk} when $q\ge 2$.

\par

If instead $1\le q\le 2$, then the assumption $r\ge d/q'$ together
with the same arguments as above show that \eqref{WFfjfk} holds
also in this case. Summing up we have
\begin{multline*}
WF_{\mathscr FL^q_{\sigma}}(\mathcal G)\subseteq \cup _{j,k=1}^N
\big (WF_{\mathscr FL^q_{\sigma}}(f_j) \cup WF_{\mathscr
FL^q_{\sigma}}(f_jf_k)\, \big )
\\[1ex]
\subseteq \cup _{j=1}^N WF_{\mathscr FL^q_{\sigma +r}}(f_j),
\end{multline*}
and (2) follows in the case $m=2$.

\par

For general $m\ge 2$, the assertion (2) now follows by repeating
these arguments and induction. This completes the proof.
\end{proof}

\par

Next we  discuss semi-linear equations. Let $J_k f $ denote
the array of all derivatives of order $ \alpha, $ $ |\alpha| \leq
k,$ of $f$ (the so called $k-jet$  of  $f$):
$$
J_k f = \{ (\partial^\alpha f) \}_{ |\alpha| \leq k}.
$$
We denote the elements of $J_k f$ with $ f_1, \dots, f_N$, where
$N$ is the number of elements in $J_k f$. We also consider $ G(x,
J_k f)$, where $G$ is the same as in Theorem \ref{F}. We also let
$P(x,D)$ to be the partial differential operator whose symbol $
P(x,\xi) $ is of the form
\begin{equation}\label{Pdef}
P(x,\xi) = \sum_{|\alpha| \leq n} b_\alpha (x) \xi^{\alpha},\quad
\text{where}\quad b_\alpha \in C^\infty (X), |\alpha| \leq n,
\end{equation}
where $X\subseteq \rr d$ is open. Note that $ P(x,D) $ is properly
supported. (Cf. Appendix A for details and notations concerning
partial differential operators, pseudo-differential operators and
sets of characteristic points.)

\par

\begin{thm}\label{semi-linear}
Let $X\subseteq \rr d$ be open, $ q \in [1,\infty] $, $r,s \geq
d/q'$, and consider the semi-linear differential equation
\begin{equation} \label{semilinearna}
P(x,D) f = G (x, J_k f),
\end{equation}
where $G$ is the same as in Theorem \ref{F}, $a _{\alpha}\in
\mathscr F L^{1} _{2s-d/q', loc} (X)$, and $P$ is given by
\eqref{Pdef}.

\par

Assume that $(x_0, \xi_0)\notin \Char (P)$, $f \in \mathscr D'
(X)$ is a solution of \eqref{semilinearna} and that one of the
following conditions hold:
\begin{enumerate}
\item[(1)] $f \in \mathscr F L^{q} _{s+k+(m-1)r, loc}(X)$,
and
\begin{equation} \label{hipo-uslov}
s+n \ge d/q' + k + (m - 1)r \text ;
\end{equation}

\vrum

\item[(2)] $f \in \mathscr F L^{1} _{s+k, loc}(X)$.
\end{enumerate}
Then $ (x_0, \xi_0) \not\in WF_{\mathscr F L^{1} _{2s+n -
d/q'}}(f)$.
\end{thm}

\par

\begin{proof} As in the proof of Theorem \ref{F},
it is enough to consider the case when $a_\alpha$ are constants.

\par

Assume that (1) holds and that we already know that $ (x_0, \xi_0)
\not\in WF_{\mathscr F L^{q} _{\sigma +k+(m-1)r}} (f) $ for some $
\sigma \geq s$. By \eqref{WFderiv} we get
$$
(x_0, \xi_0) \not\in  WF_{\mathscr F L^{q} _{\sigma + (m-1)r}} (
f^{(\alpha)} ),\quad |\alpha| \leq k,
$$
which implies that $ (x_0, \xi_0) \not \in WF_{\mathscr F L^{q} _{s + ( m-1)r}}
(f^{(\alpha)}) $, $ |\alpha| \leq k$. Hence Theorem \ref{F} (2) gives
$$
(x_0, \xi_0) \not\in  WF_{\mathscr{F} L^{q} _{\sigma}} (G (x, J_k f))
= WF_{\mathscr{F} L^{q} _{\sigma}} (P(x,D) f),
$$
provided $ \sigma \leq 2s -d/q'$. Hence, Proposition A.1 (2) in
Appendix A implies that $ (x_0, \xi_0) \not\in WF_{\mathscr F L^{q}
_{\sigma + n}} (f) $. Since $n >k$ we have gained in regularity.
Repeating the argument we obtain that $ (x_0, \xi_0)  \not\in
WF_{\mathscr F L^{q} _{2s+n - d/q'}} (f) $ as claimed.

\par

If the condition (2) holds, we repeat the above arguments for $ q
= 1$.
\end{proof}

\par

\begin{rem} Note that Theorem \ref{semi-linear}
(2) gives a hypoellipticity result which does not depend on the
order of $G$.
\end{rem}

\par

\begin{rem} In \cite[Theorem 8.4.13.]{Hrm-nonlin} it is assumed
that $G(x, J_k f)$ is in $C^\infty (\rr d)$. In that proof,
Parseval's equality and Sobolev embedding theorems are used.
Here we explain why we can not recover this case in our setting.

\par

For $2\leq q< \infty$ we have assumed that $f \in \mathscr F L^{q}
_{s}(\rr d),$ for certain $s$. Assume that $N=1$ and $G$ is given
by \eqref{velikoF}, where $a_\alpha \in C^n (\rr d)$ are such that
$\partial ^{\beta}a_\alpha$ are bounded for each $|\alpha |\le n$,
$n>d/q'$, and that $f\in H^{q'}_{(n)} (\rr d)$ instead. Here
$H^{q}_{(n)} (\rr d)$ is the Sobolev space which consists of all
$f\in L^q(\rr d)$ such that all derivatives of $f$ up to the order $n$
belong to $L^q$. Then we can show that $ G(x,f(x)) \in \mathscr F
L^{q} _{n}(\rr d).$ Namely, by \cite[Corollary 6.4.5]{Hrm-nonlin},
those assumptions imply that $\partial^{\alpha} G(f) \in L^{q'} $
for every $\alpha$, $ | \alpha | \leq m$. Now, the Hausdorff-Young
inequality gives
$$
\| \xi^{\alpha} \mathscr F(G(\cdo ,f)) \|_{L^q } = \| \mathscr F
(\partial^{\alpha} G(\cdo ,f)) \|_{L^q} \leq C \|
\partial^{\alpha} G(\cdo ,f) \|_{ L^{q'}} < \infty,
$$
for any $ |\alpha | \leq m$, hence  $ G (f) \in \mathscr F L^{q}
_{m}(\rr d)$.

\par

But still we are not able to prove Theorem \ref{semi-linear}  if $
q > 1 $ and $ G \in C^\infty (\rr d)$ (for $ G = G(f_1, \dots,
f_N)$). The reason is that we have to make localizations and
consider different constants $C$ which can not be controlled (we
do not have Taylor expansion) and, apart from this, we are not
able to control the number of appearances of $r> d/q' $.

\par

However, if $ q  =1$ and if we assume that $ G $ is real analytic,
then Theorem \ref{semi-linear} (2) can be improved. Namely,
\end{rem}

\par

\begin{prop} \label{Fprop}
Let $X$ be an open set in $\rr d$,  $s\geq 0$, $ G(y_1, \dots, y_N)$
be a real analytic function, $ f_j \in  \mathscr F L^{1} _{s, loc}
(X) $, $ j = 1,\dots, N$, and let $\mathcal G$ be the same as in
\eqref{gpoldef}. Then the following is true:
\begin{enumerate}
\item[(1)]
$\mathcal G\in \mathscr F L^{1} _{s, loc} (X) $;

\vrum

\item[(2)] If $ (x_0 , \xi_0 ) \not\in WF_{ \mathscr F L^{1} _{\sigma}} (f_j) $, $
j = 1,\dots, N$, $ s \leq \sigma \leq 2s$, then $ (x_0 , \xi_0
)\not\in WF_{\mathscr F L^{1} _{\sigma}} (\mathcal G) $;

\vrum

\item[(3)] Let $k$ and $n$ be the same as in Theorem
\ref{semi-linear}, $(x_0 , \xi_0 )\not\in WF_{\mathscr F L^{1}
_{s+k}} (f)$ where  $f$ is a solution of \eqref{semilinearna}. If
$P$ is noncharacteristic at $ (x_0, \xi_0)$ then it follows that
$(x_0, \xi_0) \not\in WF_{\mathscr F L^{1} _{2s+n}} (f) $.
\end{enumerate}
\end{prop}

\par

\begin{proof}
(1) For
$$
\mathcal G = G (f_1, \dots, f_N) = \sum_{\alpha \in
{\bf N}^N} a_{\alpha} f_1 ^{\alpha_1} \cdots  f_N ^{\alpha_N},
$$
we apply
$$
(\alpha_1 -1)+ \dots + (\alpha_N -1) + N -1
$$
times multiplication
in $ f_1 ^{\alpha_1}, \dots,   f_N ^{\alpha_N} $. By the inspection
of corresponding proofs of Theorem \ref{product} and Theorem \ref{F}
for $ q =1 $, after the joint localization for all $
f_1, \dots,   f_N  $,  we see that
$$
\int_{\rr d} | \mathscr F(\fy G (f_1, \dots, f_N) ) (\xi) |
\langle \xi \rangle^{s} d\xi \leq \sum_{\alpha \in {\bf N}^N}
|a_{\alpha}| C^{|\alpha| - 1} < \infty ,
$$
for some constant $C>0$ which depends
on $f_1,...,f_N$. This proves (1).

(2) Now we use Theorem \ref{microlocal-functions}, Theorem \ref{F}
and their proofs to conclude that there exists a  constant $ C $
depending now on localizations of $ f_1, \dots,   f_N  $ in $ x $
and $ \xi $ variables (choose cones $ \Gamma $ and $ \Gamma_1 $
adopted to any $ f_1, \dots,   f_N  $) so that 
\begin{multline*}
\int_{\Gamma} | \mathscr F(\fy G (f_1, \dots, f_N) ) (\xi) |
\langle \xi \rangle^{s} d\xi
\\[1ex]
\leq \sum_{\alpha \in {\bf N}^N} |a_{\alpha}| \int_{\Gamma} |
\mathscr F(\fy G(f_1 ^{\alpha_1}, \dots, f_N ^{\alpha_N})) (\xi) |
\langle \xi \rangle^{s} d\xi
\\[1ex]
\leq \sum_{\alpha \in {\bf N}^N} |a_{\alpha}| C^{|\alpha| - 1} < \infty.
\end{multline*}
This completes the proof of (2).

\par

Finally, the assertion (3) follows by (2).
\end{proof}

\par

\section{Modulation spaces} \label{sec6}

\par

In this section we restate our results in terms of modulation spaces.
We start with the definition of the short-time Fourier transform.

\par

Assume that $\fy \in \mathscr S'(\rr d)$ is fixed. Then the
short-time Fourier transform of $f\in \mathscr S'(\rr d)$ with
respect to $\fy$ is defined by
$$
(V_\fy f)(x,\xi ) =\mathscr F(f\cdot \overline {\fy (\cdo -x)})(\xi ).
$$
Here the left-hand side makes sense, since it is the partial
Fourier transform of  tempered distribution
$ F(x,y)=(f\otimes \overline \fy)(y,y-x) $
with respect to the $y$-variable. We also note that if $f,\fy \in
\mathscr S(\rr d)$, then $V_\fy f$ takes the form
\begin{equation}\label{stftformula}
V_\fy f(x,\xi ) = (2\pi )^{-d/2}\int f(y)\overline {\fy
(y-x)}e^{-i\scal y\xi}\, dy.
\end{equation}

\par

Assume that $\omega \in \mathscr P(\rr {2d})$, $p,q\in
[1,\infty ]$, and that $\fy \in \mathscr S(\rr d)\setminus 0$. Then
the \emph{modulation space} $M^{p,q}_{(\omega )}(\rr
d)$ consists  of all $f\in \mathscr S'(\rr d)$ such that
\begin{equation}\label{modnorm}
\begin{aligned}
\nm f{M^{p,q}_{(\omega )}} &= \nm f{M^{p,q,\fy }_{(\omega )}}
\\[1ex]
&\equiv \Big ( \int \Big ( \int |V_\fy f(x,\xi )\omega (x,\xi )|^p\,
dx \Big )^{q/p} \, d\xi \Big )^{1/q} <\infty
\end{aligned}
\end{equation}
(with obvious interpretation when $p=\infty$ or $q=\infty$). The space
$M^{p,q}_{(\omega )}$ is a Banach space which is independent of the
choice of $\fy \in \mathscr S(\rr d) \setminus 0$, and different $\fy$ give rise
to equivalent norms, \cite{F1}. Furthermore, $M^{p,q}_{(\omega )}$ increases with
respect to the parameters $p$ and $q$, and decreases with $\omega$, in
the sense
\begin{multline}\label{embmod}
M^{p_1,q_1}_{(\omega _1)}(\rr d) \subseteq M^{p_2,q_2}_{(\omega
_2)}(\rr d),\quad \text{and}\quad \nm f{M^{p_2,q_2}_{(\omega
_2)}}\le \nm f{M^{p_1,q_1}_{(\omega _1)}},
\\[1ex]
\text{when}\ p_1 \le p_2,\ q_1\le q_2,\ \omega _2\le C\omega
_1\ \text{and}\ f\in \mathscr S'(\rr d).
\end{multline}
We refer to \cite{F1,Feichtinger3,Feichtinger4,Feichtinger5,Gro-book}
and the references therein for basic facts about modulation spaces.

\par

Locally, the spaces $\mathscr FL^q_{(\omega)} (\rr d)$ and
$M^{p,q}_{(\omega )}(\rr d)$ coincide, in the sense that
$$
\mathscr FL^q_{(\omega )} (\rr d)\cap \mathscr E '(\rr d) =
M^{p,q}_{(\omega )}(\rr d) \cap \mathscr E'(\rr d),
$$
and
\begin{equation}\label{FLmodcompest}
C^{-1}\nm f{\mathscr FL^q_{(\omega )}}\le \nm f{M^{p,q}_{(\omega
)}}\le C\nm f{\mathscr FL^q_{(\omega )}},\quad f\in \mathscr E'(\rr d),
\end{equation}
for some positive constant $C$, which only depends on $d$ and the
size of the support of $f$ (see Theorem 2.1 and Remark 4.4 in
\cite{RSTT}). This property is extended in \cite{PTT} in the context
of the new type of wave front sets. In particular, for any open set
$X \subset \rr d$, since $\mathscr FL^q_{(\omega )}(\rr d)\subseteq
\mathscr FL^q_{(\omega),loc}(X)$, we have that
\begin{equation}\label{FLmodlocemb}
M^{p,q}_{(\omega )}(\rr d) \subseteq \mathscr FL^{q}_{(\omega
),loc}(X).
\end{equation}
Note that  we may recover (\ref{incrFLloc}) as a consequence of
\eqref{embmod} and \eqref{FLmodcompest}.

\par

We will also use the following (global) embeddings
\begin{multline} \label{potapanja}
M^{p_1,q}_{(\omega )} (\rr d) \subseteq \mathscr FL^q_{(\omega
)}(\rr d)  \subseteq M^{p_2,q}_{(\omega )}(\rr d),\quad \text{when}
\\[1ex]
\ p_1 \leq \min (q, q')\ \text{and}\ \max (q, q') \leq p_2.
\end{multline}
Here $\omega (x,\xi )\in \mathscr P(\rr {2d})$ is constant with
respect to the variable $x$. (Cf. e.{\,}g. Proposition 1.7 in
\cite{Toft2} and Theorem 3.2 in \cite{Toft35}.)

\par

It is convenient to set $ M^{p,q}_{s,t} (\rr d) = M^{p,q}_{(\omega )} (\rr d) $
when $\omega (x,\xi )=\eabs x^t\eabs \xi ^s$,
and we omit the index of the weight when $\omega \equiv 1$, i.{\,}e. we set
$ M^{p,q} (\rr d)  =M^{p,q}_{(\omega )} (\rr d) $ when $
\omega \equiv 1$.

\par

We remark that characterization of the wave front set
$ WF_{\mathscr FL^q_{(\omega )}}(f)$ in terms of modulation spaces and of Wiener
amalgam spaces is given in \cite{PTT}.
In particular, for $\omega \in \mathscr P(\rr {2d})$ such that
$\omega _0(\xi )=\omega (y_0,\xi )$ for some $y_0\in \rr d$ we have
\begin{equation} \label{istifrontovi}
WF_{\mathscr FL^q_{(\omega _0)}}(f)\equiv WF_{\mathscr FL^q_{(\omega
)}}(f) \equiv WF_{M^{p,q}_{(\omega )}}(f), \quad  f\in \mathscr D'(\rr
d).
\end{equation}

\par

Therefore, results of previous sections can be restated in terms of
modulation spaces. For example, the following result is an immediate
consequence of Theorem \ref{microlocal-functions},
\eqref{FLmodlocemb} and \eqref{istifrontovi}.

\par

\begin{prop} \label{micro-func-modul}
Assume that $p,q \in [1,\infty]$, and that $f_j \in  M^{\infty,q}
_{s_j,0} (\rr d)$ for $j = 1,2$. Then the following is true:
\begin{enumerate}
\item if $s_1 - |s_2| \ge 0$ when $q=1$ and $  s_1 - |s_2| > d/q'$ otherwise, then
$$
WF_{M^{p,q} _{s_2,0}}(f_1 f_2) \subseteq WF_{M^{p,q} _{|s_2|,0}}
(f_1)\text ;
$$

\vrum

\item  if instead $s_1 + s_2 \geq s \geq 0$ when $q=1$ and
$s_1 + s_2 - d/q' > s \geq 0$ otherwise, and $ s_2 - s \ge d/q' $,
then
$$
WF_{M^{p,q} _{s,0}}(f_1 f_2) \subseteq WF_{M^{p,q} _{|s_2|,0}}
(f_1)\text ;
$$
\end{enumerate}
\end{prop}

\par

Next, we consider hypoellipticity in modulation spaces. In that
context Proposition A.1 (1) and (2) are reformulated as follows.

\par

\begin{prop}\label{modul-2}
Let $ p,q \in [1,\infty]$ and  let $f \in \mathscr D'(X)$.
Then the following is true:
\begin{enumerate}
\item[(1)]  $WF_{M^{p,q} _{s-m,0}} (Af) \subset WF_{M^{p,q} _{s,0}} (f)$,
for every properly supported $ A \in \Psi^m (X)$;

\vrum

\item[(2)]
If $(x_0,\xi_0) \not\in WF_{M^{p,q} _{s-m}} (Af)$ for some
properly supported $ A \in \Psi^m (X)$ which is  noncharacteristic
at $(x_0,\xi_0)$, then $(x_0,\xi_0) \not\in WF_{M^{p,q} _{s}}
(f)$.
\end{enumerate}
\end{prop}

\par

The hypoellipticity result from Theorem  \ref{semi-linear} can be
interpreted in the following way.

\par

\begin{prop}\label{modul-semi}
Let $X\subseteq \rr d$ be open, $p,q \in [1,\infty] $ and $r,s \geq
d/q'$. Consider the semi-linear differential equation
\eqref{semilinearna} where $G$ is the same as in Theorem \ref{F}, $a
_{\alpha}\in M^{\infty,1} _{2s-d/q',0}  (\rr d)$, and $P$ is given
by \eqref{Pdef}.

\par

Assume that $(x_0, \xi_0)\notin \Char (P)$, $f \in \mathscr D' (X)$
is a solution of \eqref{semilinearna} and that one of the following
conditions are fulfilled:
\begin{enumerate}
\item[(1)] $f \in M^{\infty,q} _{s+k+(m-1)r,0} (\rr d)$,
and \eqref{hipo-uslov} hold;

\vrum

\item[(2)] $f \in M^{\infty,q} _{s+k,0}(\rr d)$.
\end{enumerate}
Then $ (x_0, \xi_0) \not\in WF_{M^{p,q} _{2s+n - d/q', 0} }(f)$.
\end{prop}

\par

Finally, for $ q =1 $ Proposition \ref{Fprop} gives the following
proposition. Here we note that $M^{\infty,1} _{s,0} (\rr d)$ is an
algebra under multiplication in view of \cite{Feichtinger0,
Feichtinger3} when $s\ge 0$.

\par

\begin{prop} \label{modul-analytic}
Let $X$ be an open set in $\rr d$,  $s\geq 0$, $ G(y_1, \dots, y_N)$
be a real analytic function, $ f_j \in M^{\infty,1} _{s,0} (\rr d)
$, $ j = 1,\dots, N$, and let $\mathcal G$ be the same as in
\eqref{gpoldef}. Then the following is true:
\begin{enumerate}
\item[(1)]
$\mathcal G\in M^{\infty,1} _{s,0} (\rr d)$;

\vrum

\item[(2)] If $ (x_0 , \xi_0 ) \not\in WF_{ M^{p,1} _{\sigma,0}} (f_j) $, $
j = 1,\dots, N$, $ s \leq \sigma \leq 2s$, then $ (x_0 , \xi_0
)\not\in WF_{M^{p,1} _{\sigma,0}} (\mathcal G) $;

\vrum

\item[(3)] Let $k$ and $n$ be the same as in Theorem
\ref{semi-linear}, $(x_0 , \xi_0 )\not\in WF_{M^{p,1} _{s+k,0}} (f)$
where  $f$ is a solution of \eqref{semilinearna}. If $P$ is
noncharacteristic at $ (x_0, \xi_0)$ then it follows that $(x_0,
\xi_0) \not\in WF_{M^{p,1} _{2s+n,0}}(f) $.
\end{enumerate}
\end{prop}

\par

\section*{Appendix A}

\par

For a detailed study of pseudo-differential operators in the
context of Fourier Lebesgue spaces we refer to \cite{PTT}. Here we
observe only the localized version of pseudo-differential
operators which is used in the study of semi-linear equations.

\par

Assume that $m\in \mathbf R$. Then we recall that the H\"ormander symbol class
$$
S^m _{1,0}= S^m_{1,0}(\rr d\times \rr d)=S^m(\rr{2d})
$$
consists of all smooth functions $a$ such that
for each pair of  multi-indices $ \alpha, \beta $ there are constants
$C_{\alpha ,\beta}$ such that
$$
| \partial^{\alpha} _{\xi}  \partial^{\beta} _{x} a (x,\xi) | \leq
C_{\alpha, \beta} \eabs \xi ^{m - |\alpha|}, \quad  x,\xi \in \rr d.
$$
We also set $S^{-\infty}_{1,0} = \cap_{m\in {\mathbf R}} S^m_{1,0}$, and
$$
\operatorname {Op}(S^m_{1,0})=\sets {a(x,D)}{a \in S^m_{1,0}(\rr d\times \rr d)},
$$
where the pseudo-differential operator $a(x,D)$ is defined by the Kohn-Nirenberg
representation
\begin{equation*}
(a(x,D)f)(x )
=
(2\pi) ^{-d} \int_{\rr d} \int_{\rr d} a(x+y,\xi )f(y)e^{i\scal {x-y}\xi }\,
dy\, d\xi.
\end{equation*}
We say that $a$ is the symbol of the operator $a(x,D)$.

\par

The symbol $a\in S^m_{1,0}(\rr {2d})$ is
called {\em non-characteristic } at $ (x_0, \xi_0) \in
\rr d \times ( \rr d \setminus 0)$
if there is a neighborhood $U$ of $x_0$, a conical neighborhood
$\Gamma$ of $\xi _0$ and constants $c$ and $R$ such that
\begin{equation}\tag*{(A.1)}
|a(x,\xi) | > c |\xi|^m, \quad \mbox{if} \quad |\xi| > R,
\end{equation}
and $\xi \in \Gamma$. Then one can find $ b \in S^{-m}_{1,0}(\rr {2d})$ such that
$$
a(x,D) b(x,D) - Id\in \operatorname{Op}(S^{-\infty}_{1,0})\quad
\text{and}\quad b(x,D)a(x,D) - Id \in \operatorname{Op}(S^{-\infty}_{1,0})
$$
in a conical neighborhood of $ (x_0, \xi_0)$ (cf. \cite{Ho1,PTT}).
The point $ (x_0, \xi_0) \in \rr d \times ( \rr d \setminus 0) $ is
called characteristic for $a$ if it is not non-characteristic point
of $a(x,D)$. The set of characteristic points   (the characteristic
set) of $a(x,D)$ is denoted by $ \Char (a(x,D))$. We shall identify
operators with their symbols when discussing characteristic sets.

\par

The operator $a(x,D) \in \operatorname{Op}(S^m_{1,0})$ is called {\em
elliptic} if the set of characteristic points is empty. This means
that for each bounded neighborhood $U$ of $x_0$, there are
constants $c,R >0 $ such that (A.1) holds when $x\in U$.

\par

A continuous linear map $ A : C_0 ^\infty (X) \rightarrow C^\infty (X) $
is said to be a pseudo-differential operator of order $m$ in $X,$ $ A \in \Psi^m (X),$
if for arbitrary $ \phi, \psi \in  C_0 ^\infty (X) $
the operator $ f \mapsto \phi A (\psi u) $ is in Op $S^m_{1,0}.$
For example, the restriction of $ a(x,D) \in $ Op$S^m_{1,0}$ to $X$ belongs to $ \Psi^m (X).$

\par

According to \cite[Proposition 18.1.22]{Ho1}, every $ A \in \Psi^m
(X)$  can be decomposed as $ A = A_0 + A_1$
where $ A_1  \in \Psi^m (X)$  is properly supported and the kernel of
$ A_0 $ is in $ C^\infty$.
In that sense it is no essential restriction to require proper
supports in the following statements.

\par

The proof of the following Theorem is based on the proof of \cite[Theorem 8.4.8]{Hrm-nonlin}
and can be found in \cite{PTT3}.

\par

\renewcommand{\rubrik}{Proposition A.1}

\begin{tom}
Assume that $ q \in [1,\infty],$ $s \in \mathbf{R}$,
$f\in\mathscr{D}'(\rr d)$ and $(x_0,\xi_0)\in \rr d \times (\rr d
\setminus 0)$. Then the following is true:
\begin{enumerate}
\item[{\rm{(1)}}] $WF_{\mathscr F L^{q} _{s - m}} (A f) \subset WF_{\mathscr{F}L^{q} _{s}} (f) $
for every properly supported $ A \in \Psi^m (X)$;

\vrum

\item[{\rm{(2)}}] if $ (x_0, \xi_0)\not\in WF_{\mathscr F L^{q} _{s - m}}( Af)$
for some properly supported $ A \in \Psi^m (X)$ which is
non-characteristic at $(x_0,\xi_0)$, then $  (x_0, \xi_0)\not\in
WF_{\mathscr F L^{q} _{s}} (f)$;

\vrum

\item[{\rm{(3)}}] there is  a conical neighborhood $U \times \Gamma $  of $( x_0,
\xi_0) $ such that $ (x, \xi) \not\in WF_{\mathscr{F} L^{q} _s} (f)$
for every $ (x, \xi) \in U \times \Gamma $ and for every $ s \in
\mathbf{R}$ if and only if $ (x_0,\xi_0) \not \in WF (f )$.
\end{enumerate}
\end{tom}

\par

As an immediate consequence of Proposition  A.1 (1) and (2) we
obtain the following.

\par

\renewcommand{\rubrik}{Proposition A.2}

\begin{tom}
Let  $ A\in \Psi^m (X)$ be properly supported. Then we have the microlocal property
$$ WF_{\mathscr{F}L^{q} _s} (f) \subset WF_{\mathscr{F}L^{q} _s} (Af) \cup \Char(A),$$
where $\Char (A)$ denotes the set of characteristic points of $A$.
\end{tom}

\par

We refer to \cite{PTT} for a more general statements of the above type.

\par

\end{document}